\documentclass[11]{article}
\usepackage{amsmath}
\usepackage{amsfonts}
\usepackage{stmaryrd}
\usepackage{graphicx}
\usepackage{cite}
\usepackage{algorithmic}
\usepackage{algorithm}
\usepackage{nicefrac}
\usepackage{exscale,relsize}

\usepackage{fullpage}

\usepackage{authblk}

\newcommand{\qed}{\,$\square$}
\newcommand{\cA}{{\mathcal A}}
\newcommand{\cN}{\mathcal{N}}
\newcommand{\cR}{{\mathcal R}}

\newcommand{\cT}{{\mathcal{T}}}

\newcommand{\cRG}{{{\mathcal R}_\Gamma}}
\newcommand{\cRGtr}{{\mathcal R}^t_{\Gamma,\psi}}
\newcommand{\cRGt}{{\mathcal R}^t_{\Gamma}}
\newcommand{\cRGst}{{\mathcal R}^*_{\Gamma}}

\newcommand{\R}{\mathbb{R}}

\newcommand{\imf}{\Gamma_\varphi^{-1}}
\newcommand{\mf}{\Gamma_\varphi}
\newcommand{\cF}{{\mathcal F}}

\newcommand{\mi}{\mathrm{i}}
\newcommand{\md}{\mathrm{d}}
\newcommand{\thetaphi}{\theta(\varphi)}

\newcommand{\cCG}{{\mathcal C}_\Gamma}

\newcommand{\cCGt}{{\mathcal C}_{\Gamma}^t}
\newcommand{\WF}{\textup{WF}}

\newcommand{\IC}{\text{IC}}

\newcommand{\cE}{\mathcal{E}}
\newcommand{\cD}{\mathcal{D}}

\newcommand{\cRGc}{{\mathcal R}_\Gamma}

\newcommand{\cQ}{\mathcal{Q}}
\newcommand{\cV}{\mathcal{V}}

\newcommand{\sigo}{\sigma_0}

\newcommand{\rhoo}{\rho_0}
\newcommand{\tx}{\tilde{x}}

\newcommand{\cl}{\operatorname{cl}}

\newcommand{\VA}{\mathcal{V}_A}
\newcommand{\IA}{\mathcal{I}_A}

\newcommand{\cotpr}{\chi_{[0,2\pi]\times \rr}}
\newcommand{\cRGr}{\mathcal{R}_{\Gamma,[0,2\pi]}}
\newcommand{\cLr}{\mathcal{L}_{[0,2\pi]}}

\newcommand{\lo}{\lambda_0}
\newcommand{\Cxo}{\mathcal{C}_{\Gamma,(x_0,\xi_0)}}
\newcommand{\cP}{\mathcal{P}}
\newcommand{\Mphi}{M_\phi}
\newcommand{\inv}{^{-1}}
\newcommand{\vpo}{\varphi_0}
\newcommand{\etao}{\eta_0}
\newcommand{\so}{s_0}
\newcommand{\xo}{x_0}

\newcommand{\Sigp}{\Sigma_{\Phi}}
\newcommand{\nn}{\mathbb{N}}
\newcommand{\st}{\,\big|\,}
\newcommand{\smo}{\setminus \boldsymbol{0}}
\newcommand{\paren}[1]{\left(#1\right)}
\newcommand{\sparen}[1]{\left\{#1\right\}}
\newcommand{\bparen}[1]{\left[#1\right]}
\newcommand{\abs}[1]{\left|#1\right|}
\newcommand{\norm}[1]{\left\|#1\right\|}

\newcommand{\bd}{\operatorname{bd}}
\newcommand{\xio}{{\xi_0}}
\newcommand{\eps}{\epsilon}

\newcommand{\vp}{\varphi}
\newcommand{\otp}{[0,2\pi]}
\newcommand{\otpr}{[0,2\pi]\times \mathbb{R}}
\newcommand{\netpe}{(-\eps,2\pi+\eps)}
\newcommand{\netper}{(-\eps,2\pi+\eps)\times \mathbb{R}}
\newcommand{\tp}{$2\pi$-}
\newcommand{\bel}[1]{\begin{equation}\label{#1}}
\newcommand{\be}{\begin{equation}}
\newcommand{\ee}{\end{equation}}
\newcommand{\cL}{{\mathcal{L}}}
\newcommand{\Gp}{\Gamma_\varphi}

\newcommand{\dx}{\mathrm{d}x}

\newcommand{\ds}{\mathrm{d}s}
\newcommand{\dsi}{\mathrm{d}\sigma}

\newcommand{\rr}{\mathbb{R}}
\newcommand{\rtwo}{\mathbb{R}^2}
\newcommand{\rn}{\mathbb{R}^n}

\newcommand{\psido}{$\Psi$DO}

\newtheorem{theorem}{Theorem}[section]
\newtheorem{definition}[theorem]{Definition}
\newtheorem{corollary}[theorem]{Corollary}
\newtheorem{proposition}[theorem]{Proposition}
\newtheorem{lemma}[theorem]{Lemma}
\newtheorem{Remark}[theorem]{Remark}
\newtheorem{Hyp}[theorem]{Hypothesis}
\newtheorem{Example}[theorem]{Example}

\begin{document}
\title{Detectable singularities from dynamic Radon data} 

\author{B.~N.\ Hahn\footnote{\textit{E-mail:}
\texttt{Hahn@num.uni-sb.de}}\ }
\author{E.~T.\ Quinto\footnote{\textit{E-mail:} 
\texttt{Todd.Quinto@tufts.edu}}}

\affil{$^*$\,Department of Applied Mathematics,
Saarland University, Saarbr\"ucken, Germany}

\affil{$^\dagger$\,Department of
 Mathematics, Tufts University, Medford, MA, USA
}


\maketitle

\begin{abstract}
In this paper, we use microlocal analysis to understand what X-ray
tomographic data acquisition does to singularities of an object which
changes during the measuring process. Depending on the motion model,
we study which singularities are detected by the measured data.  In
particular, this analysis shows that, due to the dynamic behavior, not
all singularities might be detected, even if the radiation source
performs a complete turn around the object. Thus, they cannot be
expected to be (stably) visible in any reconstruction. On the other
hand, singularities could be added (or masked) as well.  To understand
this precisely, we provide a characterization of visible and added
singularities by analyzing the microlocal properties of the forward
and reconstruction operators.  We illustrate the characterization
using numerical examples. 
\end{abstract}

\section{Introduction}\label{sect:intro}

The data collection in X-ray computerized tomography takes a certain
amount of time due to the time-dependent rotation of the radiation
source around the specimen. A crucial assumption in the classical
mathematical theory (including modeling, analysis and derivation of
reconstruction algorithms) is that the investigated object does not
change during this time period. However, this assumption is violated
in many applications, e.g. in medical imaging due to internal organ
motion. In this case, the measured data suffer from inconsistencies.
Especially, the application of standard reconstruction techniques
leads to motion artifacts in the resulting images
\cite{shepp_tuning_fork,motionartifacts}.

Analytic reconstruction methods to compensate for these
inconsistencies have been developed for special types of motion,
including affine deformations, see e.g.
\cite{crawford,roux_desbat,desbat_time_dependent}. An inversion
formula for the dynamic forward operator in case of affine motion has
been stated in \cite{Hahnaffine}, which also serves as basis for
suitable reconstruction methods. For general, non-affine deformations,
no inversion formula is known so far. Besides iterative methods, e.g.
\cite{iterative2, iterative1}, approximate inversion formulas that
accurately reconstruct singularities exist for fan beam and parallel
beam data in the plane \cite{katsevich_approx-2010} and for cone beam
data in space \cite{katsevich_local-2011}. They are based on the
observation, that operators of the form \bel{def:L-general} {\mathcal L} =
\mathcal{R}_\Gamma^t {\mathcal P} \cRG f\ee with forward operator $\cRG$,
specially designed pseudodifferential operator ${\mathcal P}$ and
backprojection operator $\cRGt$ (which is, typically, related to the
formal dual to $\cRG$), are known to reconstruct singularities of the
object. In addition, methods developed in the general context of
dynamic inverse problems have been successfully applied in
computerized tomography \cite{schmitt2,Hahnnonlinear}.

Nevertheless, there can still arise artifacts in the reconstructions,
even if the motion is known and the compensation method is exact, as
e.g. \cite{Hahnaffine}. On the other hand, the dynamic behavior of the
object can lead to a limited data problem even if the radiation source
rotates completely around the object. This means that some
singularities of the object might not be visible in the
reconstruction. 

To guarantee reliable diagnostics in practice, it is essential to
study these limitations carefully. Therefore, our aim is to analyse
which singularities are detected by the measured data in the dynamic
case and to characterize which of them can be reliably reconstructed
or whether they create additional artifacts in the reconstruction
process.

In this research, we understand the motion problem using generalized
Radon transforms and microlocal analysis. The mathematical model of
X-ray tomography with stationary specimen is integration along
straight lines \cite{Natterer86}. If the object moves during the data
acquisition, the measured data can be interpreted as data for a
(static) reference object where the integration now takes place along
curves rather than straight lines
\cite{katsevich_approx-2010,iterative2,Hahnaffine}.  Microlocal
analysis is the rigorous theory of singularities and the study of how
Fourier Integral Operators (FIO) transform them.  Guillemin
\cite{Gu1975} was first to make the connection between microlocal
analysis and Radon transforms (see also \cite{GS1977, Gu1985}) when he
showed that many generalized Radon transforms, $R$, are FIO.  He
showed that, under the Bolker Assumption (Def.\ \ref{Def:Bolker}) and
an extra smoothness assumption related to our definition of smoothly
periodic (see Sect.\ \ref{sect:periodic}), $R^*R$ is an elliptic
pseudodifferential operator (\psido). This means that $R^*R$ images
all singularities of functions and does not add artifacts.  This
theorem was exploited in \cite{Be} to show that a broad range of Radon
transforms on surfaces in $\rn$ can be ``inverted'' modulo lower order
terms.  Greenleaf and Uhlmann \cite{GU1989} and others developed the
microlocal analysis of generalized Radon transforms that occur in
X-ray CT \cite{Q1993sing, KLM}, cone beam CT \cite{FLU, Ks2006,
katsevich_local-2011}, seismics \cite{Dehoop-InsideOut}, sonar
\cite{QRS2011}, radar \cite{NC2002}, and other applications in
tomography.

Microlocal analysis has begun to be used in motion compensated CT.  In
\cite{katsevich-motion-R2-2008}, Katsevich proved that under certain
completeness conditions on the motion model, the reconstruction
operator $\cL$ in \eqref{def:L-general} detects all singularities of
the object.  This is related to theorems of Beylkin \cite{Be} showing
that operators like $\cL$ are elliptic pseudodifferential operators.
In \cite{FSU} uniqueness is proven for a broad range of Radon
transforms on curves.  The cone beam CT case is more subtle since
artifacts can be added to backprojection reconstructions, even with
stationary objects \cite{GU1989, FLU}.  Katsevich characterized the
added artifacts for this case and developed reconstruction algorithms
to, at least locally, decrease the effect of those added artifacts.
He uses this information to develop motion estimation algorithms when
the motion model is not known \cite{katsevich_local-2011}.

Motivated by large field of view electron microscopy, the article
\cite{quinto_curves} presents the microlocal analysis of general
curvilinear Radon transform as well as local reconstruction methods.
Analyzing added artifacts for X-ray tomography without motion has been
done in \cite{Katsevich:1997, FrikelQuinto2013, LN2015-limangle}, and
generalizations to other types of tomography have been done in
\cite{FrikelQuinto2015, LN2015-spherical-lim-data}.

In this article, we consider general motion models with less
restrictive completeness assumptions. To develop our characterization
of detectable and added singularities, we describe in Section
\ref{sect:math basis} the mathematical model for the dynamic case as
generalized Radon transform.  We also present the mathematical bases
of our work, including microlocal analysis.  In Section
\ref{sect:forward operator}, we assume the model is exact and study
which object singularities are encoded in the measured data.  In
Section \ref{sect:smoothly periodic} we consider the reconstruction
operator in the case of smoothly periodic motion, so the object is in
the same state at the end of the scan as the start.  Based on these
results, in Section \ref{sect:nonperiodic} we analyze the case when
limited data arise and characterize visible and added singularities in
reconstruction methods of type filtered backprojection. Our
theoretical results are evaluated on numerical examples in Section
\ref{sect:numerics}.  The more intricate proofs are in the appendix
and we show in \ref{sect:arbitrary measure} that our theorems are true
even if the weights are arbitrary on the Radon transforms.

\section{Mathematical basis}\label{sect:math basis}
We use the following notation for function spaces.  The
space of all smooth (i.e., $C^\infty$) functions of compact support is
denoted $\cD(\rn)$.  A \emph{distribution} is an element of the dual
space $\cD'(\rn)$ with the weak-* topology and pointwise convergence
(i.e., $u_k\to u$ in $\cD'(\rn)$ if for every $f\in \cD(\rn)$,
$u_k(f)\to u(f)$ in $\rr$).  Further, $\cE(\R^n)$ will denote the set
of smooth functions on $\rn$; its dual space, $\cE'(\rn)$ is the set
of distributions that have compact support.  See \cite{Rudin:FA} for a
description of the topologies and properties of these spaces.

A data set in computerized tomography can be interpreted as a function
(or distribution) with domain $[0,2\pi]\times \R$. In the static case,
the data are $2\pi$-periodic in the first variable, but this does not
necessarily hold in the dynamic case since the object does not
necessarily return to its initial state at the end of the scanning.  

Generally, smooth functions (and hence distributions) are defined on
open sets because derivatives will then be well defined. With this in
mind, we make the following definition.

\begin{definition}\label{def:smooth closed}  Let $g$ be a function with domain $[0,2\pi]\times \R \times Y$, where $Y$ is an open subset of $\R^n$.\begin{enumerate}
\item[i)] We
call $g$ \emph{smoothly periodic} if $g$ extends to a smooth function
on $\rr\times\rr\times Y$ that is $2\pi$-periodic in the first
variable.
\item[ii)] In the non-periodic case, we call $g$ \emph{smooth} if,
for some $\eps>0$, $g$ extends to a
smooth function on $(-\eps,2\pi+\eps)\times \rr\times Y$.  
\end{enumerate}  
\end{definition}

If $g$ is smoothly periodic, then $g$ can be viewed as a smooth
function on the unit circle $S^1$ by identifying $0$ and $2\pi$.  We
define $\cD(\otpr)$ as the set of all smoothly periodic compactly
supported functions on $\otpr$, and $\cD'(\otpr)$ is its dual space
with the weak-* topology.  The set of smoothly periodic functions on
$\otpr$, $\cE(\otpr)$, and its dual space $\cE'(\otpr)$ are defined in
a similar way. Including the condition of $2\pi$-periodicity in these
definitions will simplify the mapping properties of the dynamic
forward operator and its dual (see Sect.\ \ref{sect:periodic}). 

In general, the object does not return to its initial state at the end
of the scanning, i.e. its motion is not $2\pi$-periodic. For this
case, we will state our theorems and definitions using open
domains with $\vp\in \netpe$ for some $\eps>0$. Finally, distributions
can be restricted to open subsets and microlocal properties that hold
on the larger set (e.g., smoothness) hold on the smaller set.  So, our
theorems are also true when mapping to distributions on $A\times \rr$
(i.e., when the data are restricted to $A\times \rr$) when $A\subset
\netpe$ is open.

In computerized tomography with stationary specimen, the given data
correspond to integrals along straight lines of the distribution $f
\in \cE'(\R^2)$ describing the x-ray attenuation coefficients of the
investigated object. Hence, the mathematical model in the 2D parallel
scanning geometry is given by the Radon line transform
\begin{equation}
\cR f (\varphi,s) = \int_{\R^2} f(x) \, \delta(s-x^T \theta(\varphi))
\, \mathrm{d}x,
\end{equation}
with $s\in \R, \, \varphi \in [0,2\pi]$, $\theta = \theta(\varphi)=(\cos
\varphi,\sin \varphi)^T$ and the $\delta$-distribution. For fixed
source and detector position $(\varphi,s)\in [0,2\pi] \times \R$, the
integration takes place over the line
\begin{equation}
l(\varphi,s) = \lbrace x \in \R^2\st x^T\theta = s \rbrace.
\end{equation}
The data acquisition in computerized tomography is time-dependent,
since the rotation of the radiation source around the object takes a
certain amount of time. The source rotation is the only time-dependent
part of the scanning procedure since, in modern CT scanners, detector
panels are used such that all detector points record simultaneously
for each fixed source position. Concerning the mathematical model,
this means that a time instance $t$ can be uniquely identified with a
source position and vice versa. In terms of the Radon transform, the
source position is given by the angle $\varphi \in [0,2\pi]$, and
there is the unique relation to a time instance $t_\varphi \in [0, 2
\pi / \phi]$ via
\[ \varphi = t_\varphi \phi\]
with $\phi$ being the rotation angle of the radiation source. Therefore, throughout the paper, we interpret $\varphi$ also as a time instance, and $[0,2\pi]$ as time interval.   

\subsection{Mathematical model for moving objects in computerized tomography}

We now derive the mathematical model for the case when the
investigated object changes during the measuring process. A dynamic
object is described by a time-dependent function $h: [0,2\pi] \times
\R^2 \rightarrow \R^2$. In the application of computerized tomography,
$h(\varphi, \cdot) \in \cE'(\R^2)$ for a fixed time $\varphi \in [0,2\pi]$
corresponds to the x-ray attenuation coefficient of the specimen at
this particular time instance. 

The dynamic behavior of the object is considered to be due to
particles which change position in a fixed coordinate system of
$\R^2$. This physical interpretation of object movement is now
incorporated into a mathematical model.

Let $f(x) := h(0,x)$ denote the state of the object at the initial
time. We call $f$ a \emph{reference function}. Please note that $f$ is a distribution since $h(0,\cdot) \in
\cE'(\R^2)$. Further, let $\Gamma : [0,2\pi] \times \R^2 \rightarrow
\R^2$ be a motion model describing the dynamic behavior of the
specimen, where $\Gamma(0,x) = x$ and $\Gamma(\varphi,x)$ denotes
which particle is at position $x$ at the time instance $\varphi$ (in
other words, $\Gamma(\varphi,x)$ is the location at time $\varphi=0$
of the particle that is at $x$ at time $\varphi$).  For fixed $\varphi
\in [0,2\pi]$, we write \bel{def:GV}\Gamma_\varphi x :=
\Gamma(\varphi,x)\ee to simplify the notation. Using this motion model
and the reference function $f$, the state of the object at time
instance $\varphi$ is given by
\begin{equation} \label{model_referencefunction} h(\varphi,x) = f(\Gamma_\varphi x).\end{equation}

\begin{Remark}  \textup{
In the model (\ref{model_referencefunction}), each particle keeps its
initial intensity over time. However, this means that the mass of the
object may no longer be conserved. If the density varies due to the
deformation, this can be taken into account by the mathematical model
\begin{equation} h(\varphi,x) = |\det D\Gamma_\varphi^{-1} x| \,
f(\Gamma_\varphi x). \label{model_reffct_2} \end{equation} In both
cases, the respective Fourier Integral Operators describing the
dynamic setting have the same phase function and hence the same
canonical relation. Thus, our results provided in this paper hold
equivalently for the mass preserving model (\ref{model_reffct_2}), see also \ref{sect:arbitrary measure}.}\\
\end{Remark}

Throughout the paper, we make the following assumptions on the motion
model $\Gamma$, which we justify by physical properties of moving
objects and imaging systems.

\begin{Hyp} \label{hyp1} 
Let $\Gamma: [0,2\pi] \times \R^2\rightarrow \R^2$ and let
$\Gamma_\varphi x=\Gamma(\vp,x)$.  Assume $\Gamma_0 x=x$.
Then, $\Gamma$ is called a \emph{motion model} and $\Gamma_\varphi$ a
\emph{motion function} if there is an $\eps>0$ such that 
\begin{enumerate} 
\item\label{smooth} $\Gamma$  extends smoothly
to $\Gamma : (-\epsilon, 2\pi+\epsilon) \times \R^2 \rightarrow \R^2$
(so $\Gamma$ is smooth by Def.\ \ref{def:smooth closed}).
\item\label{diffeo} For each $\varphi \in (-\epsilon, 2\pi + \epsilon),
\, \Gamma_\varphi : \R^{2} \rightarrow \R^2$ is a diffeomorphism. 
\end{enumerate}
A motion model is \emph{smoothly periodic} if it satisfies these
conditions for some $\eps>0$ and  $\Gamma$ is  smoothly periodic.
\end{Hyp}

\begin{Remark}
\begin{enumerate} 
\item \textup{Note that if the motion model is smoothly periodic and
satisfies \eqref{smooth} and \eqref{diffeo} of this hypothesis for some
$\eps>0$, then it does for any $\eps>0$ because $\Gamma$ is
$2\pi$-periodic in $\vp$ in this case.}

\item \textup{Under these hypotheses, the trajectory of a fixed
particle, which is the map \[
\text{tr}_x : [0,2\pi] \rightarrow \R^2, \quad \text{tr}_x (\varphi) :=
\Gamma_\varphi^{-1} \Gamma(0,x) = \Gamma_\varphi^{-1} x,\] is a smooth curve.} \\
\textup{In practical applications of computerized tomography,
only discrete data are measured. Thus, the object's motion is
ascertained for finitely discrete time instances only, which justifies
this (theoretical) assumption of smooth trajectories.}

\item \textup{Hypothesis \eqref{diffeo}\ ensures that two particles
cannot move into the same position, and no particle gets lost (or added).
The relocation is smooth because $\Gamma$ is a smooth function.} \\

\end{enumerate}
\end{Remark}

With the mathematical model of a dynamic object
(\ref{model_referencefunction}), the operator of the dynamic setting
is given by
\bel{def:R} \cRG f (\varphi,s) := \cR (f \circ \Gamma_{\varphi})
(\varphi,s) = \int_{\R^2} f(\Gamma_{\varphi} x) \, \delta
(s-x^T\theta(\varphi)) \md x .\ee

Using the change of coordinates $z:= \Gamma_\varphi x$, we obtain the
representation
\begin{align} \label{representation_cRG}
\cRG f (\varphi,s) = \int_{\R^2} f(z) \, |\det D\Gamma_\varphi^{-1}z|
\, \delta(s-(\Gamma_\varphi^{-1}z)^T\theta(\varphi)) \, \mathrm{d}z.
\end{align}
Thus, $\cRG$ integrates the respective intensity-corrected reference
function along the curve
\begin{equation} \label{integration_curve} C(\varphi,s) =
\left\lbrace x\in \R^2 \st \left(\imf x \right)^T\theta(\varphi) = s \right\rbrace.
\end{equation}
So, for each $(\vp,s)$,  $C(\vp,s) =
\Gamma_\varphi^{-1}(l(\vp,s))$.  Because $\Gp$ is a diffeomorphism,
each $C(\varphi,s)$ is a smooth simple unbounded curve, and for each
$\vp$, the curves $s\mapsto C(\varphi,s)$ for $s\in \rr$ cover the
plane and they are mutually disjoint (they foliate the plane).

\subsection{Microlocal analysis and Fourier integral
operators}\label{sect:microlocal}

In this section we will outline the basic microlocal principles used
in the article.  We refer to \cite{Ho1971, Treves:1980PSIDO,
Treves:1980vf, Hoermander03, KrishnanQuinto2014} for more details.

The key to understanding singularities and wavefront sets is the
relation between smoothness and the Fourier transform: a distribution
$f\in \cE'(\rn)$ is smooth if and only if its Fourier transform is
rapidly decreasing at infinity. However, to make the definition
invariant on manifolds (such as $\otpr$ with $0$ and $2\pi$
identified), we need to define the wavefront set as a set in the
cotangent bundle \cite{Treves:1980PSIDO}.  So, we will introduce some
notation.  

Let $x=(x_1,\dots,x_n)\in \rn$ and $\xi=(\xi_1,\dots,\xi_n)\in \rn$.  Now let $h$ be a
smooth scalar function of variables including $x\in \rtwo$ and let
$G=(g_1,g_2)$ be a function with codomain $\rtwo$, then we define
\[ \xi\dx = \xi_1\dx_1+\cdots +\xi_n\dx_n\in
  T^*_x(\rn)\]
where $T^*_x(\rn)$ is the cotangent space at $x\in \rn$,
  \[ \partial_x h = \frac{\partial h}{\partial
      x_1}\dx_1+\frac{\partial h}{\partial x_2}\dx_2, \qquad 
    D_x h = \paren{\frac{\partial h}{\partial x_1},\frac{\partial h}{\partial
        x_2}}, \qquad G\dx = g_1\dx_1+g_2\dx_2, \]
and the
other derivatives (using $D$) and differentials (using $\partial$) are
defined in a similar way; for example, $\partial_s h = \frac{\partial
  h}{\partial s}\ds$.
  
 \begin{definition}\label{definition: Wavefront Set} Let $u\in \cD'(\rn)$ and
let $(x_{0},\xi_{0})\in \rn\times \paren{\rn\smo}$.  Then $u$ is
smooth at $\xo$ in direction $\xio$ if there is a cutoff function at
$\xo$, $\psi\in \cD(\rn)$ (i.e., $\psi(\xo)\neq 0$) and an open cone
$V$ containing $\xi_{0}$ such that $\cF(\psi u)(\xi)$ is rapidly
decreasing at infinity for all $\xi\in V$.

On the other hand, if $u$ is not smooth at $\xo$ in direction $\xio$, 
then $(\xo,\xio\dx)\in \WF(u)$, the $C^{\infty}$ wavefront set of 
$u$.\end{definition}

We now define the fundamental class of operators on which our analysis
is based: Fourier integral operators.  Note that we define them only
for the special case we use.  For other applications, one would use
the definition for general spaces in \cite[Chapter
VI.2]{Treves:1980vf} or \cite{Ho1971}.

\begin{definition}[Fourier Integral Operator (FIO)]\label{def:FIO}  
 Let $\eps>0$.  Now let $a(\vp,s,x,\sigma)$ be a smooth function on
$\netper\times \rtwo\times \rr$, then $a$ is an \emph{amplitude of
order $k$} if it satisfies the following condition.  For each compact
subset $K$ in $\netper\times \rtwo$ and $M\in \nn$, there exists a
positive constant $C_{K,M}$ such that \bel{def:symbol
order}\abs{\frac{\partial^{n_1}}{\partial
\vp^{n_1}}\frac{\partial^{n_2}}{\partial s^{n_2}}
\frac{\partial^{n_3}}{\partial
x_1^{n_3}}\frac{\partial^{n_4}}{\partial x_2^{n_4}}
\frac{\partial^{m}}{\partial \sigma^{m}} a(\vp,s,x,\sigma)}\leq
C_{K,M}(1+\abs{\sigma})^{k-m}\ee for $n_1+n_2+n_3+n_4\leq M$, $m\leq
M$ and all $(\vp,s,x)\in K$ and all $\sigma\in \rr$.

The real-valued function $\Phi\in C^{\infty}\paren{\netper\times
\rtwo\times \paren{\rr\smo}}$ is called a \emph{phase function} if 
 $\Phi$ is positive homogeneous of degree $1$ in $\sigma$ and  both 
$(\partial_{(\vp,s)}\Phi,\partial_{\sigma}\Phi)$ and
$(\partial_{x}\Phi,\partial_{\sigma}\Phi)$ are nonzero for all
$(\vp,s,x,\sigma)\in \netper\times \rtwo\times\rr\smo$. 
The phase function $\Phi$ is called \emph{non-degenerate} if on the
zero-set 
\bel{def:sigma_phi}\Sigp=\sparen{(\vp,s,x,\sigma)\in \netper\times
\rtwo\times \rr\smo\st \partial_{\sigma}\Phi=0}\ee one has that
$\partial_{\vp,s,x}\paren{\frac{\partial
\Phi}{\partial \sigma}}\neq 0$.
In this case, the operator $\cT$ defined for $u\in
\cE'(\rtwo)$ by \bel{genlFIO} \cT u(\vp,s)=\int
e^{i\Phi(\vp,s,x,\sigma)}a(\vp,s,x,\sigma)u(x)\md x\, \md\sigma\ee is
a \emph{Fourier Integral Operator (FIO)} of order $k-1/2$.  The
\emph{canonical relation} for $\cT$ is \bel{def:canonical relation}
C:=\sparen{\paren{\vp,s,\partial_{(\vp,s)} \Phi(\vp,s,x,\eta); x,
-\partial_x \Phi(\vp,s,x,\sigma)}\st (\vp,s,x,\sigma)\in
\Sigma_{\Phi}}.  \ee
\end{definition}

Note that since the phase function $\Phi$ is non-degenerate, the sets
$\Sigp$ and $C$ are smooth manifolds.  Because of the conditions on
$a$ and $\Phi$, $\cT:\cD(\rtwo)\to \cE(\netper)$ and $\cT:\cE'(\rtwo)\to
\cD'(\netper)$ is continuous in both cases \cite{Treves:1980vf}.  If the
amplitude $a$ and phase function $\Phi$ are smoothly periodic, then
the conditions in this definitions are valid on $\otpr\times
\rtwo\times \rr$ where $0$ and $2\pi$ are identified.  In this case,
$\cT u$ is $2\pi$ periodic in $\vp$ for all $u\in \cE'(\rtwo)$.

To state the theorems that form the key to our proofs, we need the
following definitions.  Let $X$ and $Y$ be sets and let $B\subset
X\times Y$,  $C\subset
Y\times X$, and $D\subset X$.  Then, 
\bel{def:composition}\begin{gathered}C^t=\sparen{(x,y)\st (y,x)\in
C}\\ C\circ D = \sparen{y\in Y\st \exists x\in D, (y,x)\in C}\\
B\circ C = \sparen{(x',x)\in X\times X\st \exists y\in Y, (x',y)\in
B,\ (y,x)\in C}.
\end{gathered}\ee  We will use these rules for sets of cotangent
vectors to calculate wavefront sets.

 \begin{theorem}[{\cite[Theorem 4.2.1]{Ho1971}}]\label{theorem:Ft} Let
$\cT$ be an FIO with canonical relation $C$.  Then the formal dual
operator, $\cT^*$ to $\cT$ is an FIO with canonical relation
$C^t$.\end{theorem}

FIO transform wavefront sets in precise ways, and our next theorem, a
special case of the H\"ormander Sato Lemma, is a key to our analysis. 

\begin{theorem}[{\cite[Theorems 2.5.7 and 2.5.14]{Ho1971}}]\label{theorem:HS} 
Let $\cT$ be an FIO with canonical relation $C$.  Let $f\in \cE'(\rtwo)$.
Then $\WF(\cT f)\subset C\circ \WF(f)$.
\end{theorem}

\def\mapsw#1{ \llap{$\myvcenter{\hbox{$\scriptstyle#1$}}$}
\mathlarger{\swarrow}}
\def\mapse#1{\mathlarger{\mathbf{\searrow}}\rlap{$\myvcenter{\hbox{$\scriptstyle#1$}}$}}
\newcommand{\myvcenter}[1]{\ensuremath{\raisebox{6pt}{\hbox{#1}}}}

To understand the more subtle properties of FIO, we investigate the
mapping properties of the canonical relation $C$.  Let $\Pi_L:C\to
T^*(\netper)\smo$ and $\Pi_R:C\to T^*(\rtwo)\smo$ be the natural
projections.  Then we have the following diagram:
\bel{micro-diagram-genl} \def\mapsw#1{
\llap{$\myvcenter{\hbox{$\scriptstyle#1$}}$} \mathlarger{\swarrow}}
\def\mapse#1{\mathlarger{\mathbf{\searrow}}
\rlap{$\myvcenter{\hbox{$\scriptstyle#1$}}$}} {\begin{matrix}& C&\\
&\mathlarger{\mapsw{\mathsmaller{\Pi_L}} \
\qquad\mapse{\mathsmaller{\Pi_R}}}&\\
{T^*(\netper)\smo} &&{T^*(\rtwo)\smo}&\end{matrix}} \ee 

First, note that if $B\subset T^*(\rtwo)$ and $D\subset T^*(\netper)$
then \bel{circ and projections} C\circ B = \Pi_L\paren{\Pi_R\inv
(B)}\qquad C^t\circ D = \Pi_R\paren{\Pi_L\inv (D)}.\ee These
statements are proven using the definitions of composition and the
projections.

 \begin{definition}\label{Def:Bolker} Let $\cT$ be an FIO with canonical
relation $C$.  Then, $\cT$ satisfies the \emph{Bolker Assumption} if
the projection $\Pi_L$ is an injective immersion.
\end{definition}

 Recall that an immersion is a smooth map with injective differential.
Victor Guillemin \cite{Gu1975, GS1977} named this assumption after
Ethan Bolker who gave a similar assumption for finite Radon
transforms.  

 \begin{definition}\label{def:elliptic} The FIO $\cT$ in
\eqref{genlFIO} is elliptic of order $m-1/2$ if its amplitude, $a$, is
of order $m$ and satisfies, for each compact set
$K\subset \netper\times \rtwo$ there are constants $C_{K}>0$ and
$S_{K}>0$ such that for all $(\vp,s,x)\in K$ and $\abs{\sigma}>S_{K}$,
$\abs{a(\vp,s,x,\sigma)}\geq C_{K}(1+\abs{\sigma})^m$.
\end{definition}

Now, we apply these ideas to dynamic tomography.

\section{Microlocal analysis of the dynamic forward
operator}\label{sect:forward operator}

In this section, we study the microlocal properties of the
forward operator $\cRG$ in dynamic computerized tomography. We show
that it is an FIO and provide conditions under which it fulfills the
Bolker Assumption. Theorem 3.6 gives the relationship between
singularities of $f$ and those of $\cRG f$ which is then analysed in
more detail, especially with respect to the importance of the Bolker
Assumption. Our theorems are true for more general FIO, but the
proofs are easier in our special case.

We now introduce some notation and describe its geometric meaning.
Here $\Gamma$ is a motion model that satisfies Hyp.\ \ref{hyp1} and
let $\eps$ be as in that hypothesis.  For $x\in \R^2, \ \varphi \in
\netpe$ define
\begin{align} \label{abbreviation_H}
H(\vp,x) := \left(\imf x\right)^T\theta(\varphi) .
\end{align}
Then, the integration curve $C(\vp,s)$ in \eqref{integration_curve} can
be written 
\[C(\varphi,s) =
\left\lbrace x\in \R^2 \st H(\vp,x) = s \right\rbrace.\]
Now, define \bel{def:N}\cN(\vp,x) :=  \partial_xH(\vp,x).\ee

Our next lemma gives the geometric meaning of this covector.

\begin{lemma}\label{lemma:normal vector}
Let $(\vpo,\so)\in \netper$ and let $x$ be a point on the integration
curve $C(\vpo,\so)$.  The vector $D_xH(\vpo,x)$ is normal the curve
$C(\vpo,\so)$ at $x$, and therefore the covector $\cN(\vpo,x)$ is conormal to
this curve at $x$.
\end{lemma}

\textbf{Proof:} The curve $C(\vpo,\so)$ is defined by the equation
$g(x):=H(\vpo,x)-\so=0$.  Therefore the gradient in $x$ of $g$ at each
$x\in C(\vpo,\so)$, which is $D_xH(\vpo,x)$, is normal to this curve
at $x$.  So, its dual covector, which is $\cN(\vpo,x)$, is conormal to
$C(\vpo,\so)$ at $x$ (i.e., in the conormal space of $C(\vpo,\so)$
above $x$).  \qed

\subsection{The canonical relation of $\cRG$}\label{sect:R}

We first prove that the forward operator (\ref{representation_cRG})
for the dynamic setting is an elliptic FIO.

\begin{theorem} \label{theorem_fio} \quad
Under Hypothesis \ref{hyp1}, the operator $\cRG$ is an elliptic FIO of
order $-1/2$ with phase function
\begin{equation} \label{phase_fct_rep} \Phi(\varphi,s,x,\sigma) :=
\sigma(s-\left(\imf x \right)^T\theta(\varphi))\end{equation} and
amplitude
\begin{equation}\label{amplitude-R} a(\varphi,s, x,\sigma) := (2\pi)^{-1} \, |\det D \imf x|
 \end{equation} which is elliptic of order zero.
\end{theorem}

The proof is given in the appendix \ref{proof:fio}.\\

Since $\cRG$ is an FIO, we can determine its canonical relation using
Definition \ref{def:FIO}, eq. \eqref{def:canonical relation}.

\begin{lemma}\label{lemma:canonical relation}
Under Hypothesis \ref{hyp1}, the canonical relation of $\cRG$ is
\bel{def:C}\begin{aligned} \cCG := \Big\lbrace
\big(\varphi,H(\vp,x),\sigma \paren{\md s - \partial_\varphi
H(\vp,x)}&; x, \sigma \cN(\vp,x)\big)\st \\ & \varphi
\in (-\eps,2\pi+\eps), \, x \in \R^2, \, \sigma \in \R\setminus
\lbrace 0 \rbrace \Big\rbrace,\end{aligned}\ee  where $\eps$ is as
given in Hypothesis \ref{hyp1}.

If the motion model is smoothly periodic in $\vp$ then the condition
on $\vp$ in \eqref{def:C} is replaced by $\vp\in [0,2\pi]$ and $\cCG$
is still a smooth manifold without boundary when $[0,2\pi]$ is
identified with the unit circle, $S^1$. 
\end{lemma}

\textbf{Proof:} According to Definition
\ref{def:FIO}, \eqref{def:canonical relation}, the canonical relation
of $\cRG$ is given by \[\cCG := \left \lbrace
(\varphi,s,\partial_{(\varphi,s)} \Phi(\varphi,s,x,\sigma);
x,-\partial_x \Phi(\varphi,s,x,\sigma)) \st (\varphi,s,x,\sigma) \in
\varSigma_\Phi \right\rbrace\] where $\varSigma_\Phi:= \left\lbrace
(\varphi,s,x,\sigma) \in (-\eps,2\pi+\eps)\times \R \times \R^2 \times
\R\setminus 0 \st \partial_\sigma (\varphi,s,x,\sigma) =
0\right\rbrace$. Using the representation of the phase function
(\ref{phase_fct_rep}) along with \eqref{abbreviation_H},
$\partial_\sigma \Phi = \paren{s-H(\vp,x)}\dsi$, and thus
$(\varphi,s,x,\sigma) \in \varSigma_\Phi$ if $s=H(\vp,x)$. The
representation of $\cCG$ then follows from the representation of the
differentials $\partial_{(\varphi,s)} \Phi(\varphi,s,x,\sigma) =
-\sigma \partial_\varphi H(\vp,x) + \sigma \md s$ and $\partial_x
\Phi(\varphi,s,x,\sigma) = -\sigma \partial_x H(\vp,x)=-\sigma
\cN(\vp,x)$, as noted in the proof of Theorem \ref{theorem_fio}.\qed

In the following theorem, we find conditions on the motion model under
which $\cRG$ satisfies the Bolker Assumption. 

\begin{theorem} \label{theorem_immersion_condition}  Assume the motion
model satisfies Hypothesis \ref{hyp1}.

\begin{enumerate}

\item \label{bolker:injective} If, for each $\vp\in
(-\eps,2\pi+\eps)$, the map \begin{equation}
\label{injectivity_condition} x \mapsto
\begin{pmatrix} H(\vp,x) \\ D_\vp H(\vp,x)
\end{pmatrix}\end{equation} is one-to-one, then $\Pi_L$ is injective.

	\item \label{bolker:immersion}If the motion model fulfills the
condition
\begin{equation} \label{immersion_condition} \IC
(x,\varphi) :=\det\begin{pmatrix} D_x H(\vp,x) \\ D_x D_\varphi H(\vp,x)
\end{pmatrix}\neq 0\end{equation} for all $x\in \R^2, \, \varphi \in (-\eps,2\pi+\eps)$, then the
projection $\Pi_L : \cCG \rightarrow T^*((-\eps,2\pi+\eps) \times \R) \setminus
0$ is an immersion.  
\end{enumerate}
Thus, under these two conditions, $\cRG$ satisfies the Bolker
Assumption (Definition \ref{Def:Bolker}).

If the motion is smoothly periodic, then $(-\eps,2\pi+\eps)$ can be
replaced by $[0,2\pi]$ in this theorem.

\end{theorem}

To illustrate the geometric meaning of condition
(\ref{injectivity_condition}) for the motion model, we assume there
exist two points $x_1$ and $x_2$ with $H(\varphi,x_1) =
H(\varphi,x_2)$ and $D_\varphi H(\varphi,x_1) = D_\varphi
H(\varphi,x_2)$ for some $\varphi \in [0,2\pi]$. The first equality
implies that the two points are on the same integration curve, i.e.
the data at angle $\varphi$ cannot distinguish between them. The
second equality means, if the angle of view $\varphi$ changes
infinitesimally, also the new curve cannot distinguish the two points
because they both stay on the same curve (at least infinitesimally).
An example for a motion model not satisfying
(\ref{injectivity_condition}), is any dynamic behavior, where two
particles, which are on the same integration curve for a time instance
$\varphi$, are rotated with the same speed and in the same direction
as the radiation source

Condition (\ref{immersion_condition}), also referred to as an
\emph{immersion condition}, is equivalent to the condition
\[ D_\varphi D_x H(\vp,x) \notin \text{span}
D_x H(\vp,x) .\] The property $\IC(x_,\varphi) = 0$ means
that, at least infinitesimally at $\varphi_0$, the line normal to the
curve\hfil\newline $C(\varphi_0,H(\varphi_0,x_0))$ at $x_0$ is stationary at
$\varphi_0$, i.e. the curves near $C(x_0,H(\varphi_0,x_0))$ are
infinitesimally rigid at $x_0$ (these statements are justified in a
related case in \cite[Remarks 2 and 5]{quinto_curves}). 

We should remark that the conditions in Theorem
\ref{theorem_immersion_condition} are essentially equivalent to the
conditions of Theorem 2.1 in \cite{katsevich-motion-R2-2008} for the
fan beam case.  There is an additional assumption in his theorem that
ensures that  all singularities are visible. 

\textbf{Proof of Theorem \ref{theorem_immersion_condition}}  On the set $\cCG$, we
introduce global coordinates $(\varphi,x,\sigma)$ by the map
\bel{def:c}\begin{aligned} c :& (-\eps,2\pi+\eps) \times \R^2 \times
\R \setminus 0 \rightarrow \cCG\\ & (\varphi,x,\sigma) \mapsto
(\varphi, H(\vp,x), \sigma(- \partial_\varphi H(\vp,x) + \md
s),x,\sigma \cN(\vp,x)).\end{aligned}\ee In these coordinates, the
projection $\Pi_L$ is given by
\begin{equation} \label{Pi_L_coordinates} \Pi_L (\varphi,x,\sigma) =
(\varphi, H(\vp,x), -\sigma \partial_\varphi H(\vp,x) +
\sigma \md s).\end{equation}  

 Using the representation (\ref{Pi_L_coordinates}) of $\Pi_L$,
one sees that $\Pi_L$ is injective if for each $\varphi \in
(-\eps,2\pi+\eps)$, the map in \eqref{injectivity_condition} is
injective.

 The map $\Pi_L$ is an immersion if its differential has constant
rank $4$, and in coordinates
\[ D \Pi_L = \begin{pmatrix}
1 & 0 & 0 & 0 \\D_\varphi H(\vp,x) & D_{x_1} H(\vp,x) & D_{x_2} H(\vp,x) & 0\\ 
-\sigma \, D_\varphi D_\varphi H(\vp,x) & -\sigma D_{x_1} D_\varphi H(\vp,x) & -\sigma D_{x_2} D_\varphi H(\vp,x) & -D_\varphi H(\vp,x) \\ 0 & 0 & 0 & 1
\end{pmatrix}.\]
Thus, condition (\ref{immersion_condition}) is equivalent to $\det
D\Pi_L \neq 0$ for all $x, \, \varphi$, and thus is equivalent to,
$\Pi_L$ being an immersion.\qed

The importance of this Bolker Assumption for the detection of object
singularities in dynamic Radon data is discussed in the next section.

\subsection{Visible Singularities}\label{sect:visible singularities}

Now, we classify singularities of functions that appear in the data,
both algebraically and geometrically.

\begin{theorem} \label{theorem_wavefrontsets}
Assume the motion model, $\Gamma$, satisfies Hypothesis \ref{hyp1}.
Let $f\in \cE'(\rtwo)$.  Then, 
\begin{equation} \label{wavefrontset_general_motion} \WF (\cRG f) \subset \cCG \circ \WF (f). \end{equation}
Now assume, in addition, that
$\cRG$ satisfies the Bolker Assumption.  Then, 
\begin{equation} \label{wavefrontset_bolker} 
\WF (\cRG f) = \cCG \circ \WF (f). \end{equation}
\end{theorem}

We will prove this theorem in the appendix, \S
\ref{proof:wavefrontsets}.

The explicit correspondence between object and data
singularities is given in the following corollary.

\begin{corollary}\label{corollary_detectable_sing}
Let $f \in \cE'(\rtwo)$, and let $\Gamma$ be a motion model satisfying
Hypothesis \ref{hyp1}.  Let $A$ be an open subset of $\netpe$ and let
$(\varphi_0,s_0) \in A\times \rr$, $\sigma\neq 0$, $\beta\in \rr$.

If $(\varphi_0,s_0; \sigma (\md s - \beta \md \varphi))
\in \WF(\cRG f)$ then 
there is an $\xo\in C(\vpo,\so)$  such that  \[(x_0, \sigma 
\cN(\varphi_0,x_0) ) \in \WF(f)\] where $C(\vpo,\so)$ is the
integration curve given by \eqref{integration_curve} and $\cN$ is given
by \eqref{def:N}.

Now assume in addition
$\cRG$ satisfies the Bolker Assumption.   For $\vpo\in \netpe$,
\bel{sing correspondence}
\begin{gathered}(\varphi_0,s_0; \sigma (\md s - \beta \md \varphi))
\in \WF(\cRG f).\\
\text{if and only if}
\\
\text{there is an $\xo\in C(\vpo,\so)$  such that  $(x_0, \sigma 
\cN(\vp_0,x_0) ) \in \WF(f)$}.\end{gathered}\ee 
Furthermore, if such a point $\xo$
exists, then it is unique.
\end{corollary}

The proof follows immediately from Theorem
\ref{theorem_wavefrontsets}  and the
expression \eqref{def:C} for the canonical relation $\cCG$.
In particular, the first statement follows from
\eqref{wavefrontset_general_motion}, and the equivalence \eqref{sing
correspondence} follows from \eqref{wavefrontset_bolker}.

For $B\subset \netper$ define \bel{def:T*B}T^*_B(\netper) =
\sparen{(\vp,s,\eta)\st (\vp,s)\in B,\ \eta \in
T^*_{(\vp,s)}(\netper)}.\ee

 Corollary \ref{corollary_detectable_sing} justifies our next
definition.

\begin{definition}\label{Def:visible sing} 
Let $A\subset \netpe$ and let $\Gamma$ be a motion model satisfying
Hypothesis \ref{hyp1}.  Assume the associated Radon transform, $\cRG$,
satisfies the Bolker Assumption.  Let $f\in \cE'(\rtwo)$ and let
$(\xo,\xio)\in\WF(f)$.  Then, $(\xo,\xio)$ is a \emph{visible
singularity from data $\cRG f$ above $A$} if $\xio$ has the
representation \begin{equation}\label{representation_xi0} \xi_0 =
\sigma  \cN(\vp_0,x_0)\end{equation} for some $\sigma \neq
0$ and $\vpo \in A$.

 We call \bel{def:VA} \VA = \sparen{(x,\sigma \cN(\vp,x)\st
x\in \rtwo, \vp \in A, \sigma\neq 0}\ee the set of all possible
visible singularities from $\cRG$ above $A$. Covectors in
\[\IA=\paren{T^*(\rtwo)\smo}\setminus\cV_{\cl(A)}\]
will be called \emph{invisible singularities from $A$}. \end{definition}

Using \eqref{circ and projections}, it follows that \bel{projections
of VA} \VA = \cCGt\circ T^*_{A\times \rr}(\netper) =
\Pi_R\paren{\Pi_L\inv\paren{T^*_{A\times \rr}(\netper)}}.\ee

Corollary \ref{corollary_detectable_sing} justifies the definition: if
the motion model satisfies Hypothesis \ref{hyp1} and $\cRG$ satisfies
the Bolker Assumption, then a singularity $(x,\xi)\in \WF(f)$ causes a
singularity from the data $\cRG f$ above the open set $A$ (i.e., in
$T^*_{A\times \rr}(\netper)$) if and only if it is in $\VA$.  The
singularities of $f$ that are in $\IA$ are smoothed by $\cRG$.  Note
that the singularities of $f$ in $\cV_{\bd(A)}$ are problematic
because we will show they are in directions that can be added
singularities or that can be visible or masked by added singularities.

We can now describe the geometric meaning of the visible
singularities.

\begin{corollary}
	Let the motion model fulfill the Bolker Assumption.  The
dynamic operator $\cRG$ detects a singularity of $f$ at a point $x_0$
in direction $\xi_0$ if and only if there is an integration curve
passing through $x_0$ with  $\xi_0$  conormal to the curve at $\xo$
(i.e., the curve has tangent line at this point that is
normal to $\xi_0$).
\end{corollary}

 \textbf{Proof:} Let $\so = H(\vpo,\xo)$.  Corollary
\ref{corollary_detectable_sing} shows that, under the Bolker
Assumption, a singularity of $f$ at $(\xo,\xio)$ is visible if and
only if $\xio = \sigma \cN(\vpo,\xo)$ for some $\sigma\neq 0$.
Furthermore, Lemma \ref{lemma:normal vector} establishes that for each
$(\vp,s)\in \netper$ and each $x\in C(\vp,s)$, the covector $\cN(\vp,x)$
is conormal to $C(\vp,s)$ at $x$.  Thus a singularity of $f$ at
$(\xo,\xio)$ is visible if and only if $\xio$ is conormal to
$C(\vpo,\so)$ at $\xo$.\qed

\begin{Remark} \textup{
  In general, each data singularity at a point in data space,
$(\vpo,\so)$, stems from an object singularity $x_0\in C(\vpo,\so)$
with direction $\xi_0$, where $\xi_0$ is perpendicular to the curve
$C(\vpo,\so)$ at $\xo$.  However, in case the Bolker Assumption is not
fulfilled by the motion model, two object singularities could cancel
in the data and thus, not lead to a corresponding data singularity.\\
  In contrast, under the Bolker Assumption, every singularity in the data
  comes from a singularity in the object. Note that Example
    \ref{example_limited} shows that not all singularities of the object
    necessarily show up in the data.}
\end{Remark}

 Another way to understand visible singularities is the following.
$(\xo,\xio)\in \VA$ if there is some $\sigma\neq 0$ and $\vp_0\in A$,
such that $\xi_0 \in \text{Range}(\mu_{{x_0}})$, where $\mu_{{x_0}}$
is the map \bel{def:mu}\mu_{{x_0}}(\sigma,\varphi_0) = \sigma
\cN(\vp_0,x_0)\ee for $(\sigma,\varphi_0)\in \paren{\rr\smo}\times A$
(see (\ref{representation_xi0})).  If this map $\mu_{{x_0}}$ is not
injective, {the} object singularity ${x_0}$ can cause two different
data singularities, resulting in redundant data, as illustrated by our
next example.

 \begin{Example} \label{example_redundant} Let the dynamic behavior of
$f$ be given by the rotation $\Gamma_\varphi x= A_\varphi x$ with
rotation matrix \[ A_\varphi = \begin{pmatrix} \cos \varphi & -\sin
\varphi\\ \sin \varphi & \phantom{-} \cos \varphi
\end{pmatrix}.\]
This describes an object which rotates in the opposite direction
as the radiation source with the same rotational speed. In particular,
it holds $\Gamma_\varphi = \Gamma_{\varphi+2\pi}$ for $\varphi \in
[0,2\pi]$, so this is a smoothly periodic motion model. Since
$A_\varphi$ is a unitary matrix for all $\varphi \in [0,2\pi]$, it is
\[ H(\vp,x) = (A_\varphi^{-1}x)^T\theta(\varphi) = x^T
A_\varphi\theta(\varphi) = x^T \theta(2 \varphi).\] By a calculation
using it's definition, $\IC(x,\varphi) = 2 \cos^2 (2\varphi) + 2
\sin^2 (2\varphi) = 2$, and the map \[x \mapsto \begin{pmatrix}
x^T\theta(2\varphi) \\ 2 \, x^T\theta(2\varphi)^\perp
\end{pmatrix}\] is one-to-one since the matrix $\big(\theta(2\varphi), \theta(2\varphi)^\perp\big)^T$
is nonsingular. Thus,
the dynamic operator $\cRG$ satisfies the Bolker Assumption, and
$\WF(\cRG f) = \cCG \circ \WF(f)$.

Now, let $(x_0, \xi_0 \md x) \in \WF(f)$ with $\xi_0 := \theta(\pi)$.
Since it holds that \[ \cN(\tfrac{\pi}{2},x_0) =
\begin{pmatrix} \cos \pi \\ \sin \pi
\end{pmatrix} = \xi_0,\]
as well as \[D_x H(\tfrac{3\pi}{2},x_0) = \begin{pmatrix}
\cos \pi \\ \sin \pi
\end{pmatrix} = \xi_0,\]
this one singularity in object space causes two singularities
\begin{align*}(\pi/2,H(\pi/2,x_0),\sigma \md s - \sigma
x^T\theta^\perp(\pi)  \md
\varphi) &\in \WF(\cRG)\ \text{ and }\\ (\tfrac{3
\pi}{2},H(\tfrac{3\pi}{2},x_0),\sigma \md s - \sigma
x^T\theta^\perp(\pi) \md \varphi) &\in \WF(\cRG).\end{align*} This is according to
the fact that the projection $\Pi_R : \cCG \rightarrow
T^*(\R^2)\setminus 0$ is not injective due to the motion introduced
data redundancy. \end{Example}

If the map $\mu_{{x_0}}$ in \eqref{def:mu} is surjective for all $x_0
\in \R^2$ then \textbf{all} singularities and \textbf{all} directions
are gathered in the measured data, and we speak of \emph{complete
data}. In the static case, this corresponds to the fact that the
radiation source rotates around the complete circle (e.g.,
\cite{Q1980}). If $\mu_{{x_0}}$ is not surjective, when the point
$x_0$ is only probed by data from a limited angular range. The
following example illustrates that the dynamic behavior of the object
can lead to incomplete data, even if the full angular range $[0,2\pi]$
is covered by the source.

\begin{Example} \label{example_limited} We consider the rotational movement $\Gamma_\varphi x = A_\varphi x$ with
\[ A_\varphi = \begin{pmatrix}
\phantom{-} \cos (\tfrac{2}{3} \varphi) & \sin (\tfrac{2}{3} \varphi)\\ - \sin (\tfrac{2}{3} \varphi) & \cos (\tfrac{2}{3} \varphi)
\end{pmatrix}.\]
In this setting, the object rotates in the same direction as the
radiation source with half of its rotation speed. In particular, this
is a non-periodic motion model. It is
\[ H(\vp,x) = x^T A_\varphi\theta(\varphi) = x^T \begin{pmatrix}
\cos (\tfrac{\varphi}{3}) \\ \sin (\tfrac{\varphi}{3})
\end{pmatrix}. \]
One shows the injectivity condition, (\ref{injectivity_condition}),
is fulfilled in the same way as in Example \ref{example_redundant}.
Computing the derivatives, we obtain $\IC(x,\varphi) = \tfrac{1}{3}
\cos^2(\tfrac{\varphi}{3}) + \tfrac{1}{3} \sin^2(\tfrac{\varphi}{3}) = \tfrac{1}{3}$.  So, the Bolker Assumption
holds.

Now, assume $(x_0,\xi_0 \md x)\in \WF(f)$ with $\xi_0 =
\theta(\tfrac{5\pi}{6})$. According to Theorem
\ref{theorem_wavefrontsets}, a corresponding singularity is seen in
the data if there exists an angle $\varphi_0 \in [0,2\pi]$ with $\xi_0
= A_{\varphi_0} \theta(\varphi_0) = \theta(\tfrac{\varphi_0}{3})$, or $\xi_0 = - A_{\varphi_0} \theta(\varphi_0) = \theta(\tfrac{-\varphi_0}{3})$.
Since $\tfrac{\varphi_0}{3} \in [0,\tfrac{2}{3}\pi]$ for all $\varphi_0 \in
[0,2\pi]$, an angle $\varphi_0$ with the required property does not
exist. Hence, the singularity $(x_0,\xi_0 \md x)\in \WF(f)$ cannot be
seen in the data.\\
\end{Example}

\section{The dynamic reconstruction operator for smoothly periodic motion}
\label{sect:smoothly periodic}
In this section, we prove the main theorem  for smoothly periodic
motion. Basically, under our assumptions, the reconstruction operator is
well behaved and reconstructs all singularities of the object without
introducing new artifacts. First, we define the backprojection operator.

\subsection{Backprojection for Smoothly Periodic 
Motion}\label{sect:periodic}

In general, we denote the backprojection operator by $\cRGt$ and
define it as 
\begin{equation} \label{def:Rt} \cRGt
g(x) = \int_{\vp\in \otp} |\det D\Gamma_\varphi^{-1} x| \, g(\varphi,
(\imf x)^T\theta(\varphi)) \, \md \varphi. \end{equation} Note that, for
smoothly periodic motion, this backprojection operator is the formal
dual, $\cRGst$, to $\cRG$ for $g \in \cE([0,2\pi]\times \R)$.  A
generalization to arbitrary weights is explained in section
\ref{sect:arbitrary measure}.

\begin{proposition}\label{prop:composition}  If the motion model $\Gp$
  is smoothly periodic, then the backprojection operator, $\cRGt$,
can be composed with $\cRG$ for $f\in \cE'(\rtwo)$ and, if $\cP$ is a
pseudodifferential operator, then the reconstruction operator
\[\cL = \cRGt\cP \cRG\] is defined and continuous on domain 
$\cE'(\rtwo)$.\end{proposition}

 \textbf{Proof:} The proof will now be outlined.  First, we show when
$f\in \cD(\rtwo)$, $\cRG f\in \cD(\otpr)$.  By the smoothness
assumptions on $\Gp$, the integrals over $C(\vp,s)$ vary smoothly in
each variable, and because $\Gp$ is \tp periodic, the curves are \tp
periodic (i.e., $C(\vp+2\pi,s) = C(\vp,s)$). Thus, the integrals $\cRG
f(\vp,s)$ are smooth and \tp periodic because each $f\in \cD(\rtwo)$
has fixed compact support and $\Gp$ is \tp periodic.  Now, to show
$\cRG$ is continuous, one considers the seminorms on $\cD(\otpr)$ (see
\cite[Part II, 6.3]{Rudin:FA}).  So, assume $f_k\to f$ in
$\cD(\rtwo)$; this means that the sequence $(f_k)$ and all derivatives
converge uniformly to those of $f$, and the $f_k$ and $f$ are all
supported in a fixed compact set $K\subset \rtwo$.  By continuity of
$\Gp$ and compactness of $\otp$, there is an $R>0$ such that
$C(\vp,s)\cap K=\emptyset$ for $\abs{s}>R$, so $\cRG f_k$ and $\cRG f$
are supported in $\otp\times [-R,R]$. Finally, one uses Lebesgue's
Dominated Convergence Theorem and properties of derivatives and
integrals to show that $\cRG f_k$ and all derivatives in $\vp,s$
converge uniformly to those of $\cRG f$ and are all supported in a
fixed compact set in $\otpr$.  Since $\cRGt$ is the formal dual to
$\cRG$ in the smoothly periodic case, an analogous proof shows that
$\cRGt:\cE(\otpr)\to \cE(\rtwo)$ is continuous. 

By duality, if the motion is smoothly periodic, then
$\cRG:\cE'(\rtwo)\to \cE'(\otpr)$ and $\cRGt:\cD'(\otpr)\to \cD'(\rtwo)$
are both weakly continuous.  Since $\cP:\cE'(\otpr)\to \cD'(\otpr)$ is
also continuous, $\cL$ is weakly continuous.\qed

\subsection{The main theorem for smoothly periodic motion}\label{sect:main
  periodic}

Our main theorem for this case gives conditions under which our
reconstruction operator images all singularities and adds no
artifacts.  It is a parallel beam analogue of the fan beam result of
Katsevich \cite[Theorem 2.1]{katsevich-motion-R2-2008}. However, in
that article, the backprojection operator has a different measure; our
proof would still be valid in this case, see section
\ref{sect:arbitrary measure} of the appendix.  The same distinctions
apply to \cite{Be} and the proof outline in the last section of
\cite{KLM} for generalized Radon transforms.  Furthermore, because of
their goals, these authors consider only a few special filters, $\cP$.

\begin{theorem}\label{thm:psido-R}  Assume the motion model is
smoothly periodic and $\cRG$ satisfies the Bolker Assumption.  Let
$\cL = \cRGt \cP \cRG$ where $\cP$ is an elliptic pseudodifferential
operator with everywhere positive symbol.  Then, $\cL$ is an elliptic
pseudodifferential operator.  Therefore, for any $f\in \cE'(\rtwo)$, 
\bel{WF=}\WF(\cL F) = \WF(f).\ee \end{theorem}

The proof of Theorem \ref{thm:psido-R} will be given in  appendix
\S \ref{proof:psido-R}.

\begin{Remark}\label{smoothly periodic}  We highlight several
implications of the theorem and its proof.  

By \eqref{WF=}, all singularities are visible if the motion is smoothly
periodic and satisfies the Bolker Assumption.

Furthermore, in Remark \ref{remark:elliptic general}, we prove that
$\cL$ is elliptic as long as the pseudodifferential operator $\cP$ is
positive on $\Pi_L(C)$.  The standard Lambda tomography filter
$\cP=-d^2/ds^2$ and the standard filtered back projection operator
$\cP=\sqrt{-d^2/ds^2}$ both satisfy this condition, even though their
symbols are not elliptic on $T^*(\otpr)$.

Finally, the positivity condition can be further relaxed, and this
will be explained in Remark \eqref{remark:elliptic
general}.\end{Remark}

\section{Non-periodic motion and  added artifacts}
\label{sect:nonperiodic}  
If the motion model is smoothly
periodic and satisfies the Bolker Assumption then all singularities
are visible from the data (see Remark \ref{smoothly periodic}), and
$\cL = \cRGt \cP \cRG$ reconstructs all singularities if $\cP$ is
elliptic with positive symbol (see Theorem \ref{thm:psido-R}).
However, in smoothly periodic motion, the investigated object is in
the same state at beginning and end of the data acquisition. Thus, in
applications, this condition will in general not be met.  

In this section, we therefore study what can be said for non-periodic
motion models under the Bolker Assumption.  We assume the model
satisfies Hypothesis \ref{hyp1}, so the motion model is defined on
$\netper$ for some $\eps>0$.  However, in practice, the data are taken
only on $\otpr$. Note that the microlocal analysis developed in
Section 3 is valid on an open interval and, for non-periodic motion,
data are given on $\otpr$ This creates problems that we will now
analyze.

\subsection{The forward and backprojection operators for non-periodic
motion} Since the data are given on $\otpr$, the forward operator must
be restricted, so $\cRG$ must be multiplied by the characteristic
function of $\otpr$ to restrict to the data set.  Therefore, the
restricted forward operator is \bel{def:restricted forward} \cRGr
:=\cotpr \cRG\ee

For convenience in the proof, the backprojection operator will use the
formal dual to $\cRG$ on $\netper$ rather than $\cRGt$.  One can show
for integrable functions, $g$, that the formal dual to $\cRG$ is
defined by \begin{equation} \label{def:R*} \cRGst g(x) = \int_{\netpe}
|\det D\Gamma_\varphi^{-1} x| \, g(\varphi, (\imf x)^T\theta(\varphi))
\, \md \varphi. \end{equation} Since $\cRGst$ does not have domain
$\cD'(\netper)$, we multiply by a cutoff function.  Let
$\psi:\netpe\to \rr$ be equal to one on $\otpr$ and be supported in
$\netpe$.  We let \bel{def:restricted dual}\cRGtr g =
\cRGst\paren{\psi g}.\ee Prop.\ \ref{prop:composition nonpdic} shows
that this restricted dual is defined for $g\in \cD'(\netper)$.

 The restricted reconstruction operator is defined as \bel{def:cLr}
\cLr = \cRGtr \cP \cRGr\ee where $\cP$ is a pseudodifferential
operator in data space. In the course of the proof of Theorem \ref{thm:artifacts}
we will prove that these operators are defined for distributions and
can be composed (see Proposition \ref{prop:composition nonpdic}).

\subsection{Characterization of artifacts for the reconstruction
operator with non-periodic motion} 
In the
following, we characterize the propagation of singularities
under reconstruction in case of a non-periodic motion model. 

 Let $A\subset \netpe$, then, for $f\in \cE'(\rtwo)$, we define
\bel{def:WFA}\WF_{A}(f):= \WF(f) \cap \cV_{A} \ee where $\VA$ is
defined in \eqref{def:VA}. When $A$ is open, $\WF_A(f)$ is the set of
\emph{visible singularities of $f$} for data from $A$.  If $A$ is
closed, there can be added artifacts in the reconstruction from the
boundary $\bd(A)$, as will be shown in our next theorem.

\begin{theorem}\label{thm:artifacts}
Let $f\in \cE'(\R^2)$, and ${\mathcal P}$ be a pseudodifferential operator
and $\cLr$ is given by \eqref{def:cLr}. Then,
\begin{align*}
\WF(\cLr f) \subset \WF_{\otp}(f)\cup \cA (f),
\end{align*}
where
 \begin{align} \label{set_add_sing} \cA(f) := \lbrace & (\tx,\sigma
\cN(\tx,\varphi)): \varphi \in \lbrace 0,2\pi\rbrace,\ s\in
\rr, \ \tx\in C(\vp,s),\ \sigma\neq 0,\\
&\qquad\text{and }\   \exists \, x \in
C(\varphi,s),\ (x,\sigma (\cN(\vp,x)))
\in \WF(f) \rbrace \notag
\end{align} denotes the set of \emph{possible added artifacts}. 
\end{theorem}

\begin{Remark}\label{remark:non periodic}
This theorem shows that only singularities $(x,\xi)\in \WF(f)$ with
directions in the visible angular range can be reconstructed from
dynamic data. Singularities of $f$ outside of $\cV_{[0,2\pi]}$ are
smoothed.

Additionally, if $f$ has a singularity at a covector $(x,\sigma
\cN(\varphi_0,x))$ where $\varphi_0\in \lbrace 0,
2\pi\rbrace$, then that singularity can generate artifacts all along the
curve $C(\vpo,x)$.  These covectors are in the set 
\[\cA(\vpo,x,\sigma) = \sparen{\paren{\tx,\sigma \cN(\vpo,\tx)}\st
\tx\in C(\vpo,x)}.\] Note that the covector $\cN(\varphi_0,x)$ is
conormal to the curve $C(\vpo,x)$ at $x$ by Lemma \ref{lemma:normal
vector}.  

Furthermore, the set $\cA(f)$ is the union of the $\cA(\vpo,x,\sigma)$
for \[\vpo \in \sparen{0,2\pi}, \qquad(x,\sigma \cN(\vpo,x))\in
\WF(f).\]\end{Remark}

Under positivity conditions on $\cP$, we will also have a lower bound on the
visible singularities of $f$ that are recovered by $\cLr f$.

\begin{theorem} \label{thm:ellipticity nonperiodic}Let  $\cRG$ be a motion
  model satisfying the Bolker assumption.  Assume $\cP$ is an elliptic
pseudodifferential operator.  Finally, assume the uniqueness condition
\bel{unique covector}\begin{gathered} \forall 
(x,\xi)\in T^*(\rtwo),\ 
 \text {there is at most one $(\vp,s)\in \netper$ with}\\
\text{$x\in C(\vp,s), $ and $\xi$ conormal to $C(\vp,s)$ at
$x$}\end{gathered}\ee holds.  Then, \bel{WF= non-periodic}
\WF_{(0,2\pi)} (f)=\WF_{(0,2\pi)}(\cLr f)\ee where $\WF_{(0,2\pi)}$ is
defined in \eqref{def:WFA}.\end{theorem}

This shows that, in this case, visible singularities in
$\cV_{(0,2\pi)}$ can be recovered.  This theorem is valid under some
weaker assumptions but the statements are more technical.  The biggest
obstacle to weakening the uniqueness assumption \eqref{unique
covector} occurs when a singularity at $(x,\xi)$ is conormal to a
curve $C(\vpo,\so)$ for $\vpo\in (0,2\pi)$ \emph{and} conormal to
curves at ends of the angular range: $C(0,s_1)$ or $C(2\pi,s_2)$.
Then, added artifacts along $C(0,s_1)$ or $C(2\pi,s_2)$ could cancel a
real singularity at $(x,\xi)$.  Ellipticity theorems with more general
assumptions than \eqref{unique covector} are given for the hyperplane
transform in \cite[Theorem 5.4]{FrikelQuinto-hyperplane}.

\subsection{An artifact reduction strategy}\label{sect:smooth cutoff}

For motion that is not smoothly periodic, there is another way to
handle the limited data for $\vp$ in $\otp$ rather than multiplying by
a sharp cutoff, $\cotpr$. One can make $\cRG$ and $\cRGt$ \tp
periodic by multiplying by a smooth cutoff function, $\phi$, in $\vp$
that has compact support in $(0,2\pi)$ and is equal to one on most of
this interval.  In this case, the smoothed reconstruction operator
would be \bel{def:cLp} \cL_{\phi}(f) = \paren{\cRGt\phi} \cP
\paren{\phi\cRG f}\ee and, for $f\in \cD(\rtwo)$, $\phi\cRG f$ is
smooth and \tp periodic so in $\cD(\otpr)$.  Then, these operators can
be composed and are continuous on distributions and the proof is
essentially the same as the proof of Proposition
\ref{prop:composition}.

 Under the Bolker Assumption, $\paren{\cRGt \phi}\paren{\cP
\paren{\phi \cRG}}$ is a standard pseudodifferential operator.  The
proof is essentially the same as in the smoothly periodic case because
$\phi \cRG$ and its formal adjoint, $\cRGst\phi = \cRGt\phi$, are FIO
satisfying the Bolker assumption. 

It's important to point out that this reconstruction operator is not
necessarily elliptic everywhere, even though it is a standard
pseudodifferential operator.  Furthermore, not only the added
artifacts will be smoothed out, visible singularities near $\cA(f)$
(i.e., for covectors $(x,\eta(\vp,x)$ for $\vp$ near $0$ or $2\pi$)
will be attenuated as well because the cutoff $\phi$ is zero near $0$
and $2\pi$.  

This idea has been used in X-ray tomography without motion in
\cite{FrikelQuinto2013, FrikelQuinto-hyperplane} and generalizations
to non-smooth cutoffs are in \cite{Katsevich:1997}.  The analogous
idea is used in \cite{katsevich-motion-R2-2008} for motion compensated
CT in the fan-beam case.

\section{Numerical examples}\label{sect:numerics}

In this section, we use our theoretical results to analyze the
information content in the measured data using numerical examples.
First, we consider a specimen which performs a rotational movement
during the data acquisition, in addition to the rotation of the
radiation source, where $\Gamma_\varphi x = A_\varphi x, \ x \in \R^2,
\, \varphi \in [0,2\pi] $ with the unitary matrix from Example
\ref{example_limited} \[ A_\varphi := \begin{pmatrix} \phantom{-}\cos
(\tfrac{2}{3} \varphi) & \sin(\tfrac{2}{3} \varphi) \\ -\sin
(\tfrac{2}{3} \varphi)& \cos(\tfrac{2}{3} \varphi)
\end{pmatrix}, \quad \varphi \in [0,2\pi].\] Note that this rotation is not $2\pi$-periodic. 

\begin{figure} 
\begin{center}
\begin{minipage}[b]{.45\linewidth}
	\includegraphics[width=0.75\linewidth]{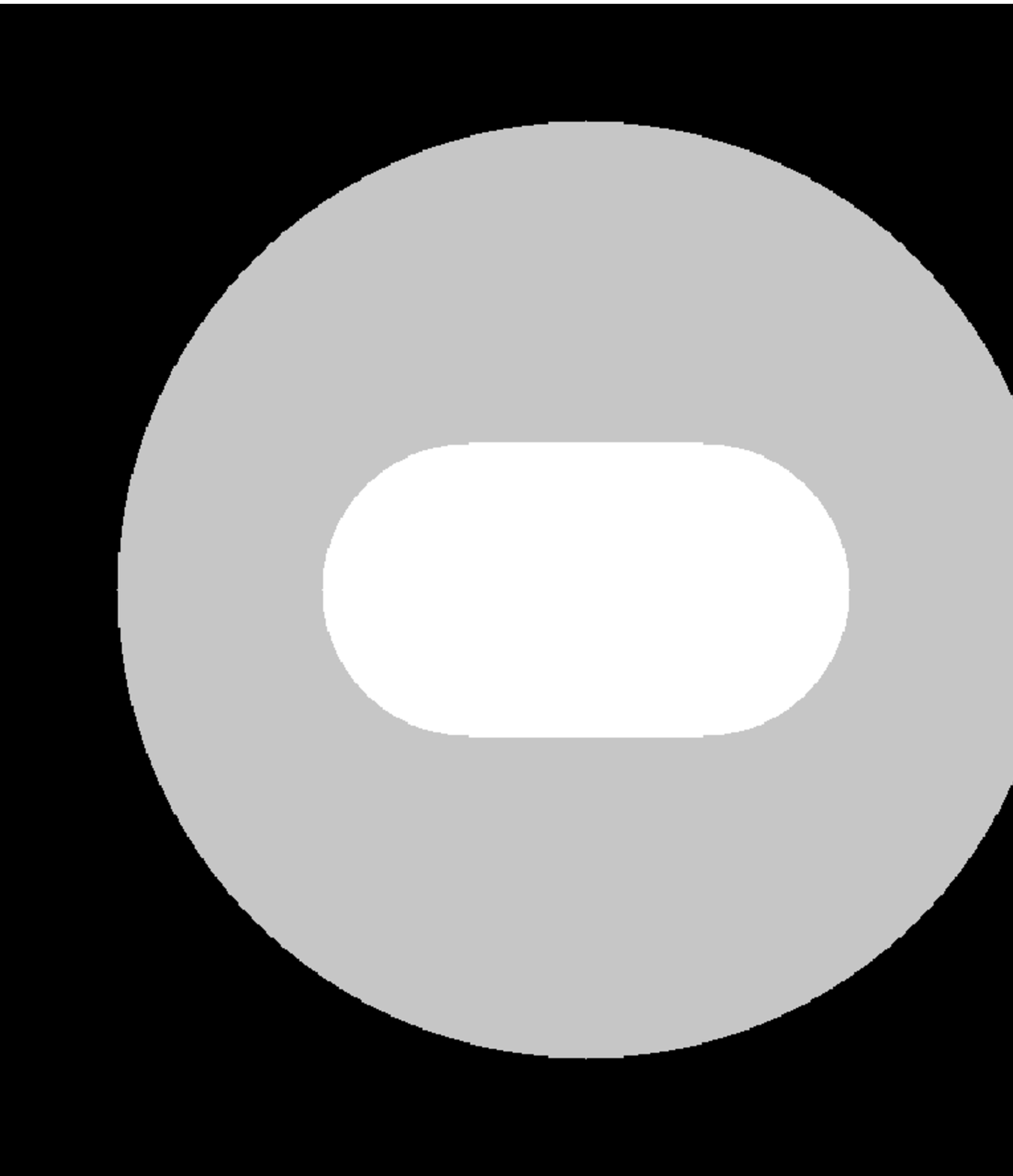}
\caption{\small Object at time instance $\varphi=0$ (reference
object)}\label{Figure_7}
\end{minipage}
\hspace{.04\linewidth}
\vspace{.005\linewidth}
\begin{minipage}[b]{.45\linewidth}
	\includegraphics[width=0.75\linewidth]{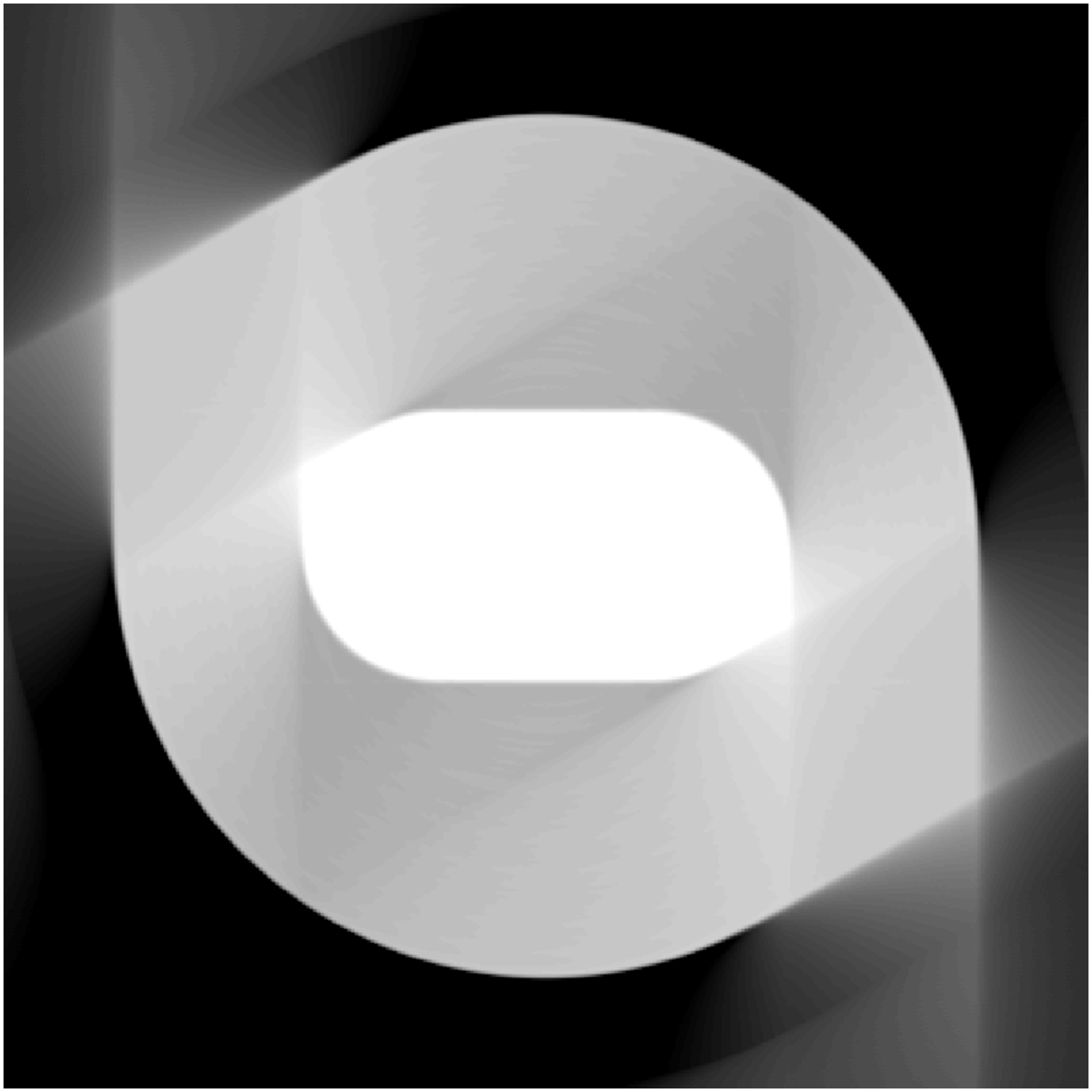}
	\caption{\small Reconstruction incorporating  exact motion functions}\label{Figure_8}
\end{minipage}
\end{center}
\end{figure}

The initial state, i.e. the reference function $f$, of our specimen is
displayed in Figure \ref{Figure_7}.  The motion corrupted Radon data
$\cRG f$ are computed in the 2D parallel scanning geometry with
$p=300$ uniformly distributed angles in $[0,2\pi]$ and $450$ detector
points.

In Example \ref{example_limited}, it is shown that not all singularities of
the specimen are ascertained by the measured data. More precisely, a
singularity $(x,\xi \, \md x) \in \WF(f)$ is detected if there is a
$\varphi \in [0,2\pi], \sigma \in \R$ such that \[ \xi_0 = \sigma
D_xH(\vp,x) = \sigma \theta(\tfrac{\varphi}{3}). \]
Thus,
\[ \lbrace \sigma D_xH(\vp,x)\st \ \varphi \in [0,2\pi], \, \sigma \in \R\setminus 0 \rbrace = \lbrace \sigma \theta(\varphi) \st \varphi \in [0,\tfrac{2 \pi}{3}] \cup [\tfrac{4\pi}{3},2\pi], \, \sigma \in \R\setminus 0 \rbrace, \]
i.e. only singularities with direction
\begin{equation}\label{ex2_detectable_sing} \xi = \sigma
\theta(\varphi_\xi), \quad \varphi_\xi \in [0,\tfrac{2}{3}\pi]\cup
[\tfrac{4\pi}{3},2\pi]\end{equation} are gathered in the data. In other
words, singularities with direction $\xi = \sigma \theta(\varphi_\xi),
\ \varphi_\xi \in (\tfrac{2}{3}\pi,\tfrac{4}{3}\pi)$
cannot be reconstructed from the dynamic data set. 

This is clearly seen in the reconstruction result, see Figure
\ref{Figure_8}. Here, we used the exact motion functions and the
algorithm proposed in \cite{Hahnaffine} as reconstruction method which
compensates known affine deformations exactly. In \cite{Hahnaffine},
it is outlined that the algorithm is of type filtered backprojection,
and hence, it fits into our framework of reconstruction operators
$\mathcal{L}_{[0,2\pi]} = \mathcal{R}^t_{\Gamma,\psi} \mathcal{P}
\mathcal{R}_{\Gamma,[0,2\pi]}$, see (\ref{def:cLr}).

Further, the singularities gathered at time instance $\varphi=0$ and $\varphi=2\pi$ create added artifacts along their integration curve. Since
\[ C(\varphi,s) = \lbrace x \in \R^2 \st (\Gamma_\varphi^{-1}x)^T\theta(\varphi) = s \rbrace = \lbrace x \in \R^2 \st x^T A_\varphi \theta(\varphi) = s \rbrace,\]
these added artifacts arise along straight lines with direction $\displaystyle{\theta\left( \tfrac{4}{3}\pi\right)^\perp}$ and $\displaystyle{\begin{pmatrix}
0\\-1\end{pmatrix}}$. Thus, the reconstructed image, Figure
\ref{Figure_8}, shows the typical limited angle streak artifacts known from the static case on the angular range $(\tfrac{2}{3}\pi,\tfrac{4}{3}\pi)$.

\begin{figure}
\begin{minipage}[b]{.24\linewidth}
    \includegraphics[width=1\linewidth]{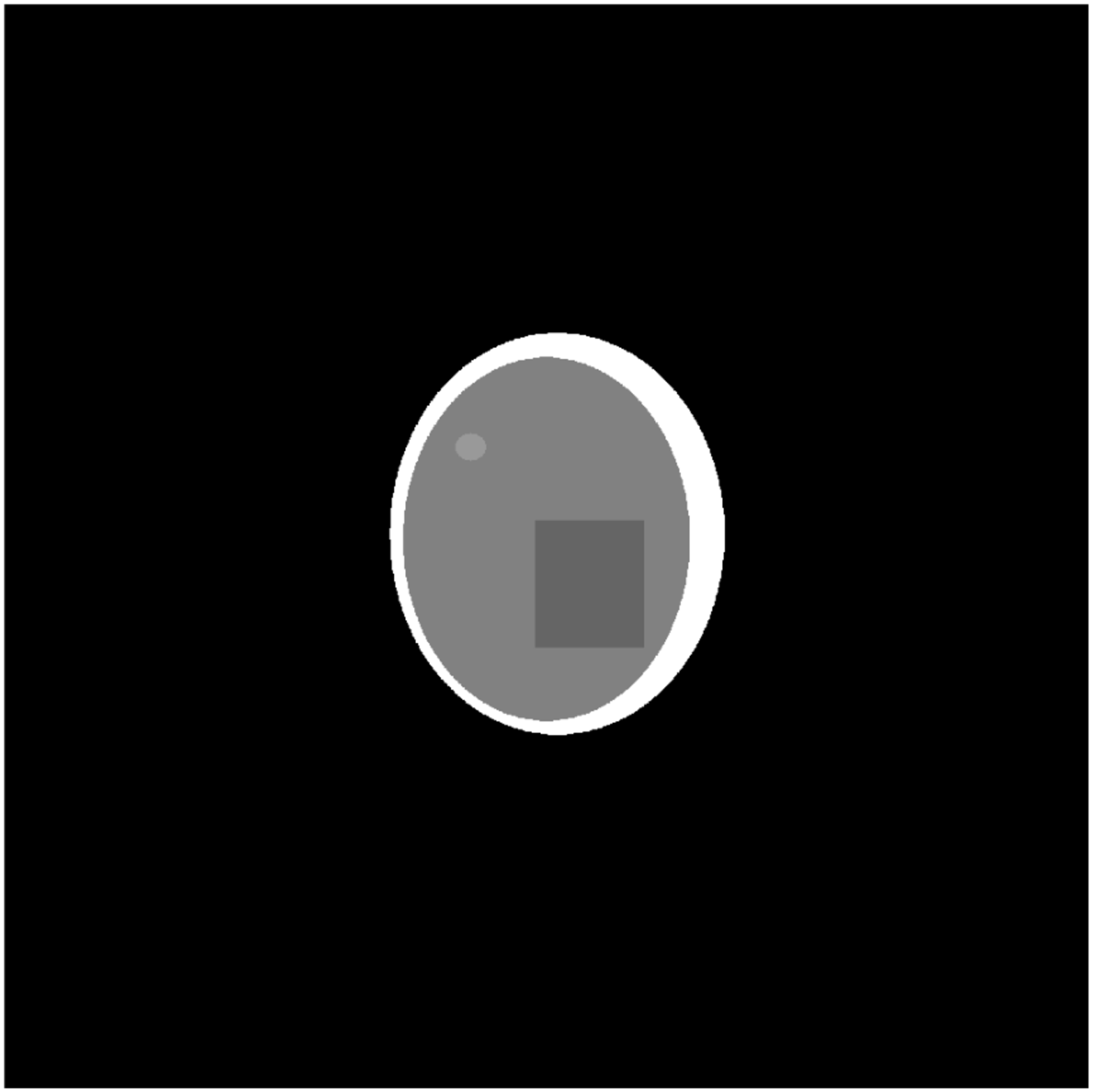} 
    \end{minipage}
    \begin{minipage}[b]{.24\linewidth}
    \includegraphics[width=1\linewidth]{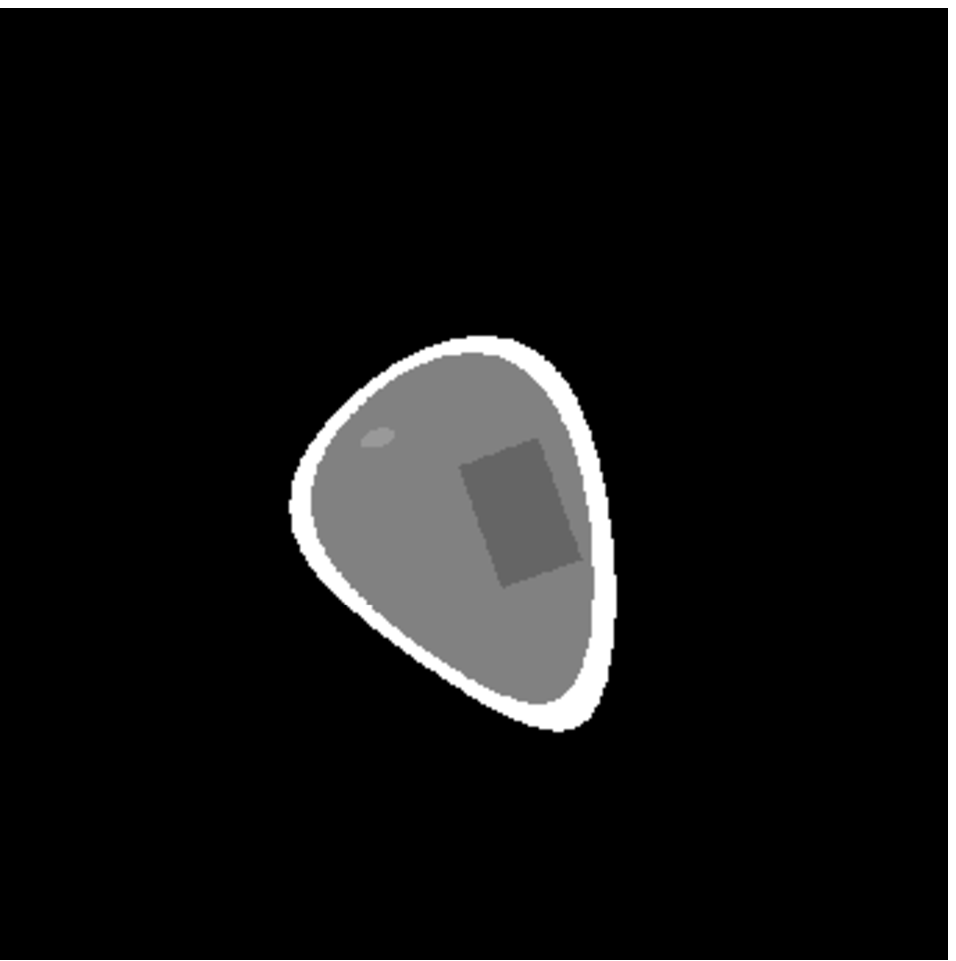} 
    \end{minipage}
    \begin{minipage}[b]{.24\linewidth}
    \includegraphics[width=1\linewidth]{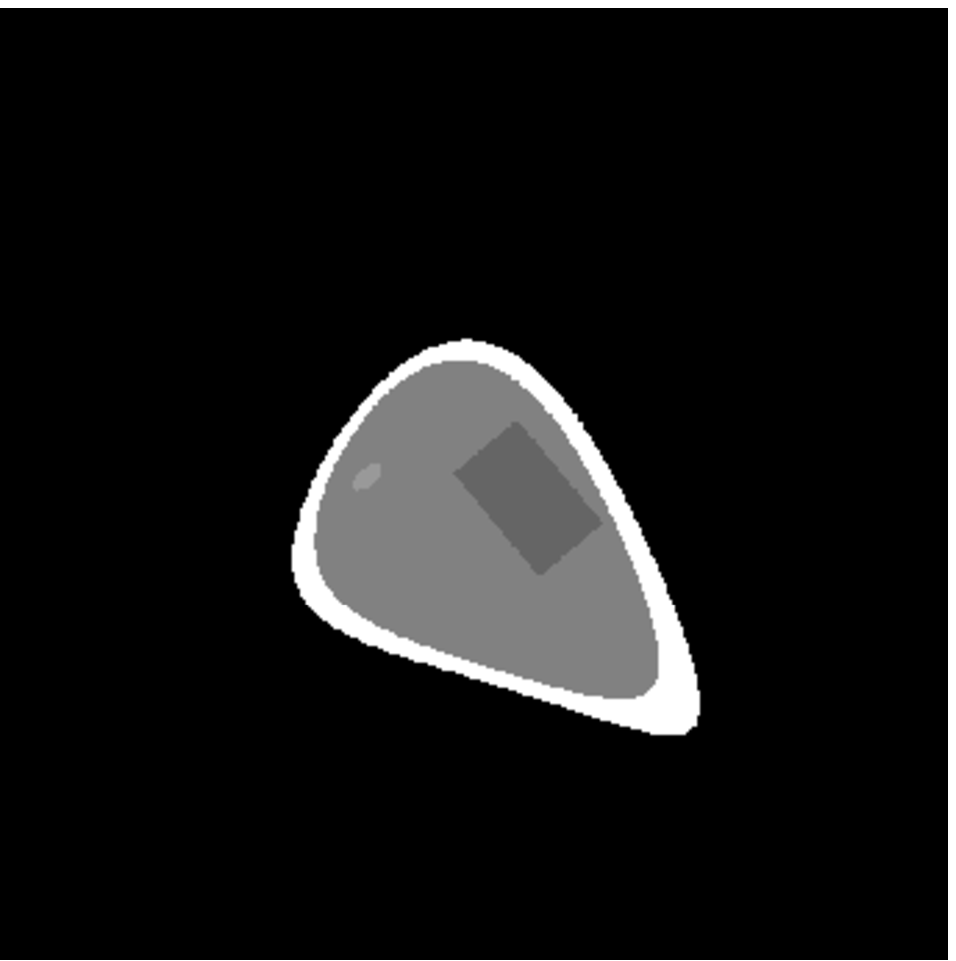} 
    \end{minipage}
    \begin{minipage}[b]{.24\linewidth}
    \includegraphics[width=1\linewidth]{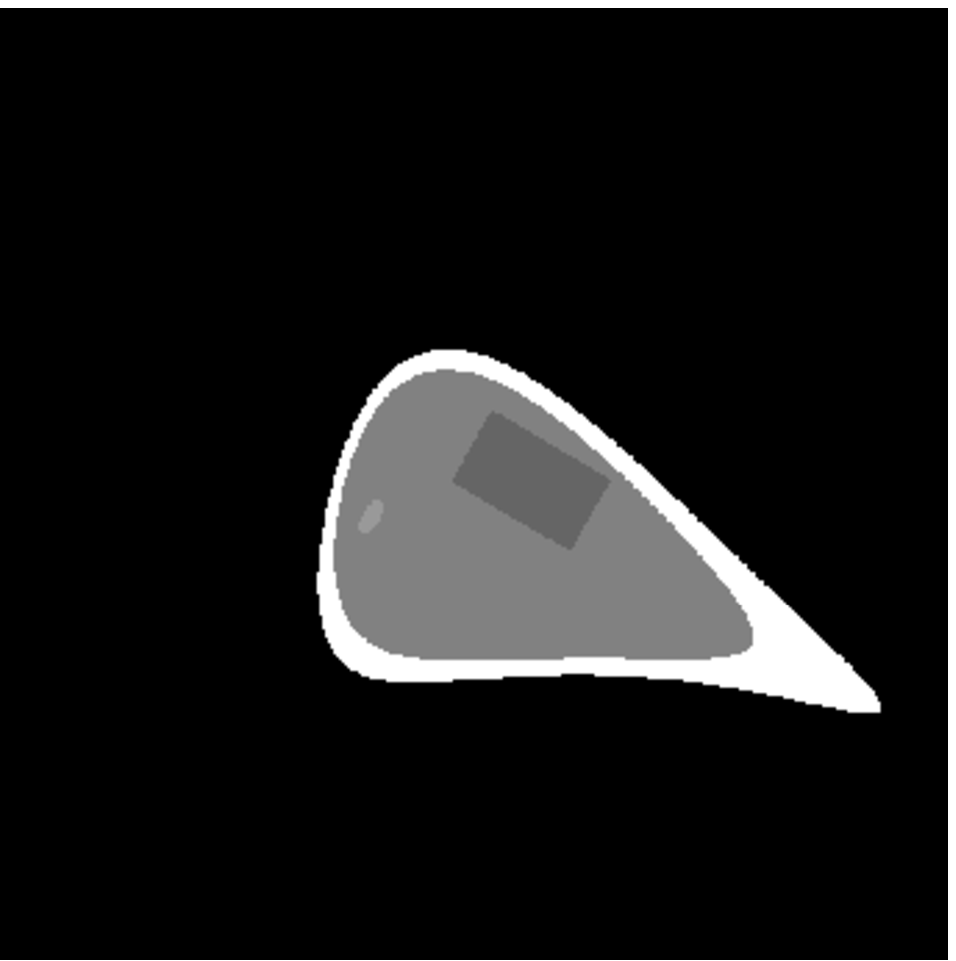} 
    \end{minipage}
    \caption{Non-affine motion of the phantom during the scanning}\label{Figure9}
    \end{figure}

Next, we illustrate our results for a non-affine motion model, where the integration curves $C(\varphi,s)$ no longer correspond to straight lines. As an example, we consider the non-periodic motion model
\[ \Gamma_\varphi x = \Gamma_\varphi^\text{scal} A_\varphi x\]
with rotation matrix
\[ A_\varphi = \begin{pmatrix}
\phantom{-}\cos (\tfrac{2}{3}\varphi) & \sin (\tfrac{2}{3}\varphi)\\ - \sin (\tfrac{2}{3}\varphi) & \cos (\tfrac{2}{3}\varphi)
\end{pmatrix}\]
and 
\[ \Gamma_{\varphi} x = \left(\begin{array}{c}
x_1 \, s_1(\varphi,x)\\x_2\, s_2(\varphi,x)
\end{array}\right)\]
with scaling parameters that depend on the time $\varphi$ as well as
on the particle $x$, see \cite{Hahnnonlinear}. In the numerical
example, \[s_i(\varphi,x) = \sum_{j=0}^4 (\sqrt[4]{5m_i} \, x_i)^j,
\quad i=1,2,\] with $m_1=\sin(5\, \cdot 10^{-5}\, \varphi \, p/\pi), \
m_2 = \sin (7\, \cdot 10^{-5} \, \varphi \, p/\pi)$.  The deformation
of the object during the data acquisition is illustrated in Figure
\ref{Figure9}.

In \cite{Hahnnonlinear}, a reconstruction method was proposed which compensates for non-affine motion, and which belongs to the class of reconstruction operators $\mathcal{L}_{[0,2\pi]} = \mathcal{R}^t_{\Gamma,\psi} \mathcal{P} \mathcal{R}_{\Gamma,[0,2\pi]}$, see (\ref{def:cLr}).

\begin{figure} 
\begin{center}
\begin{minipage}[b]{.45\linewidth}
	\begin{center}\includegraphics[width=0.75\linewidth]{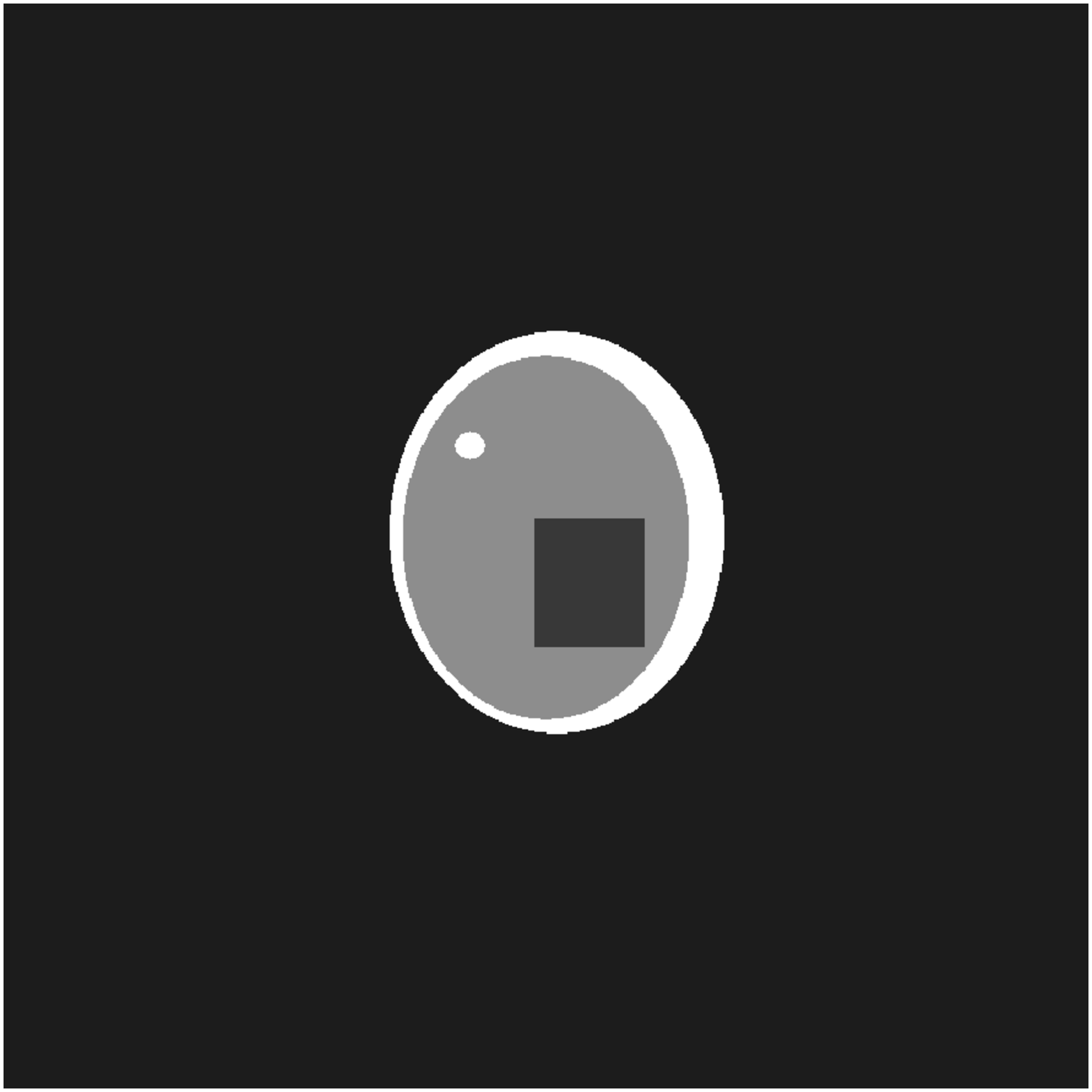}\end{center} 
\caption{\small Object at time instance $t=0$ (reference
object)}\label{Figure10}
\end{minipage}
\hspace{.04\linewidth}
\vspace{.005\linewidth}
\begin{minipage}[b]{.45\linewidth}
	\begin{center}\includegraphics[width=0.75\linewidth]{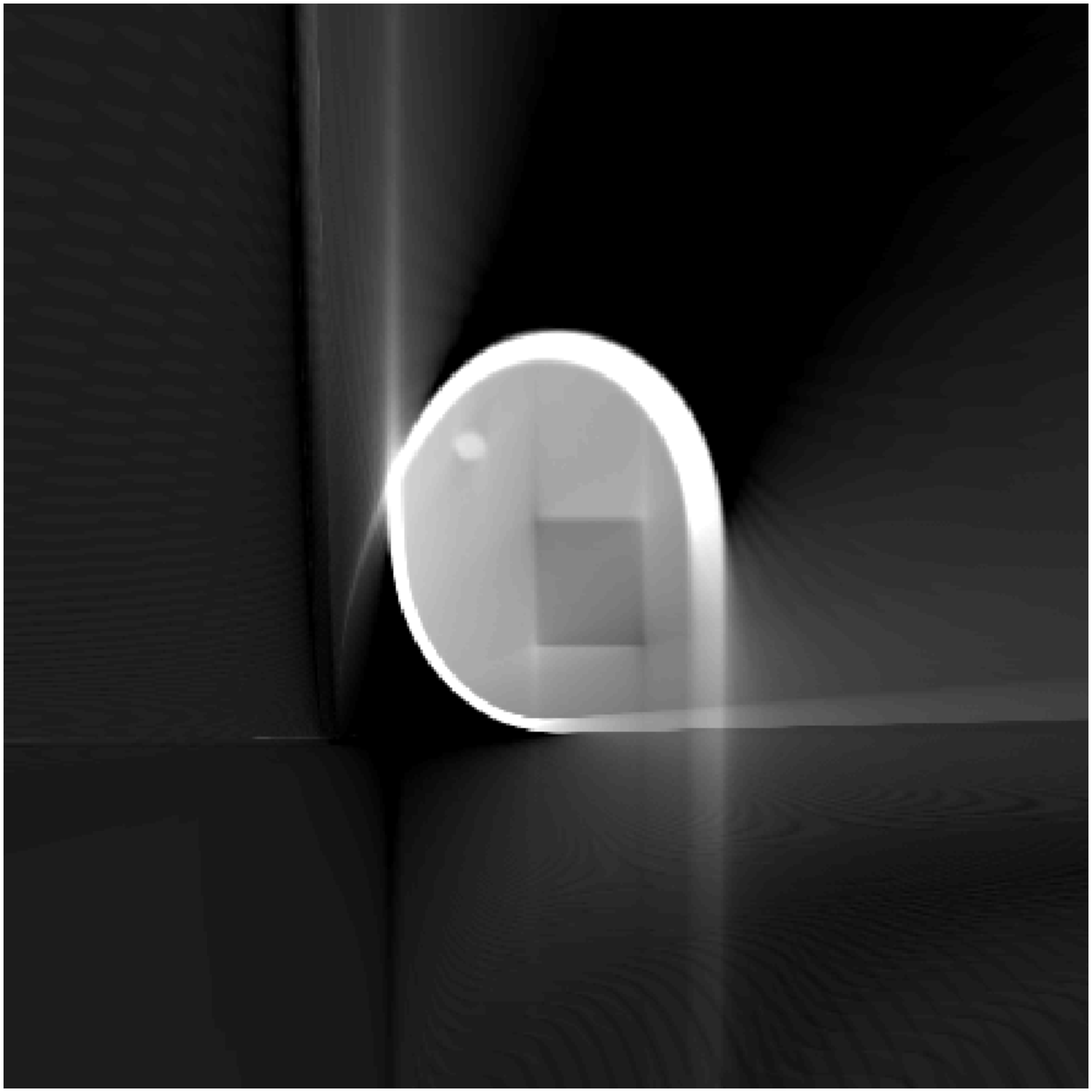} \end{center} 
	\caption{\small Reconstruction incorporating exact motion functions}\label{Figure11}
\end{minipage}
\vspace{.05\linewidth}

\begin{minipage}[b]{.45\linewidth}
	\begin{center}\includegraphics[width=0.75\linewidth]{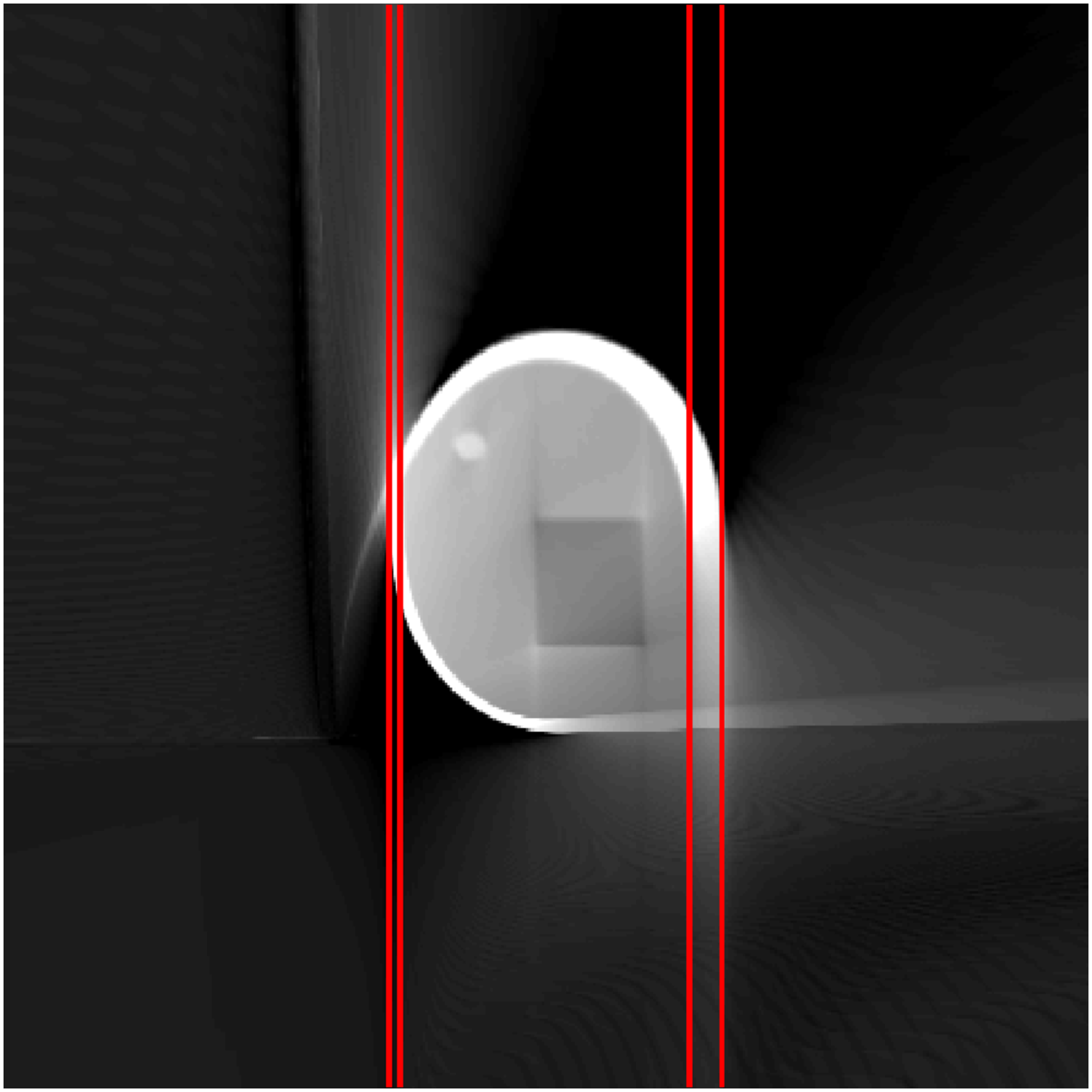}\end{center}
\caption{\small Reconstruction with integration curves at time instance $\varphi=0$ }\label{Figure12}
\end{minipage}
\hspace{.04\linewidth}
\vspace{.005\linewidth}
\begin{minipage}[b]{.45\linewidth}
	\begin{center}\includegraphics[width=0.75\linewidth]{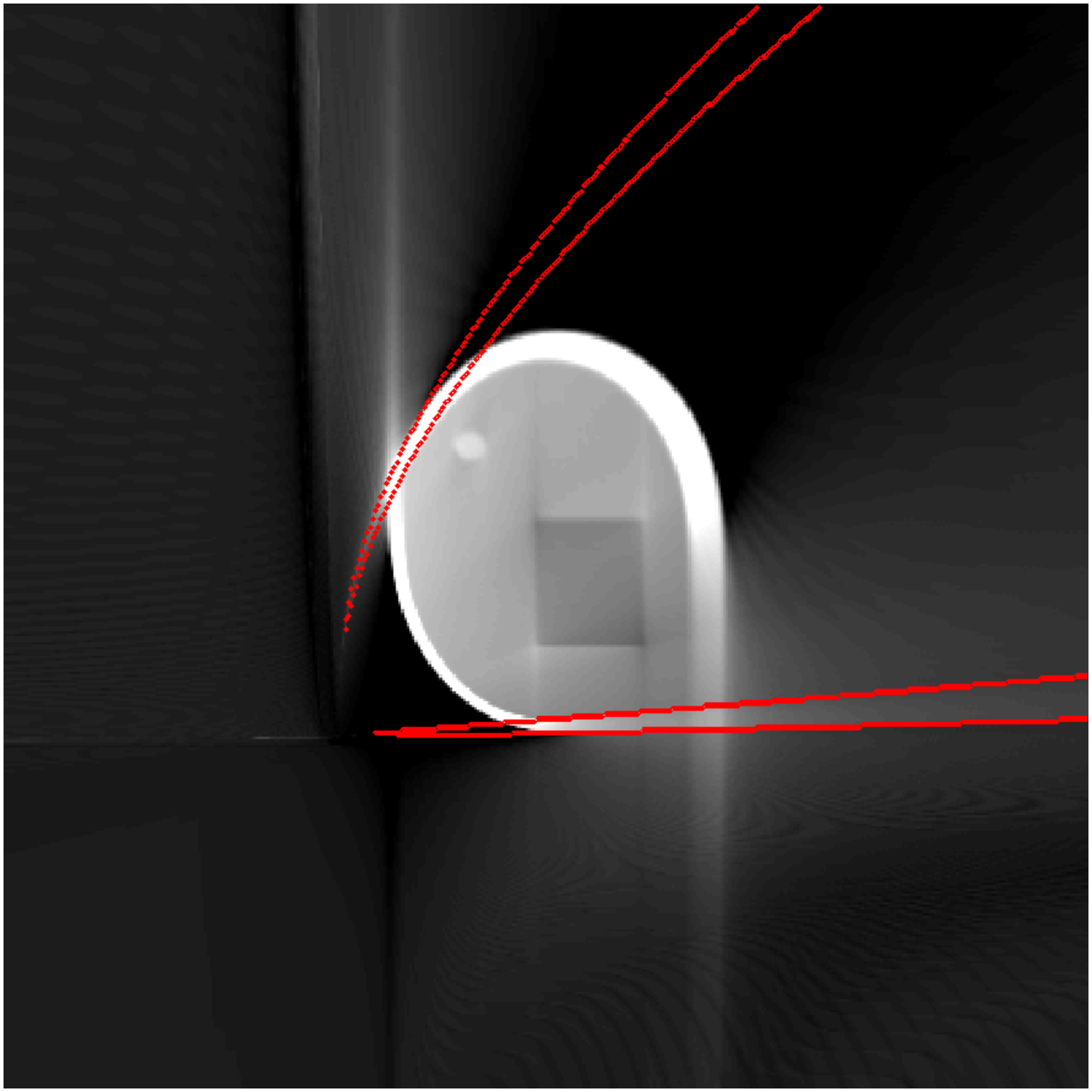}\end{center}
	\caption{\small Reconstruction with integration curves at time instance $\varphi=2\pi$}\label{Figure13}
\end{minipage}

\end{center}

\end{figure}

Applying this method to the dynamic data set provides an image showing
the visible singularities, i.e. the ones ascertained by the data, as
well as additional artifacts, see Figure \ref{Figure11}.  Figure
\ref{Figure12} and \ref{Figure13} display in addition the integration
curves passing through the singularities of the two outer ellipses,
detected at time instance $\varphi=0$ and $\varphi=2\pi$,
respectively. The comparison shows that, in accordance to our theory,
the additional artifacts spread along these integration curves. Since
$\Gamma_0 x = x$, the curves for $\varphi = 0$ are straight lines,
whereas at $\varphi = 2\pi$, they are indeed curves, not straight
lines.

\section{Conclusion and Outlook}

In this article, it was shown that the dynamic behavior of the object
in computerized tomography can lead to limited data problems, and this
means that certain singularities will be invisible in the
reconstruction results, regardless of the performance of the motion
compensation algorithm. We also provide a characterization of
detectable singularities that depends on the exact dynamics, as well
as possible added artifacts which arise even if the object's dynamic
behavior is exactly known in the reconstruction step. In applications,
this has to be taken into account at the evaluation of the
reconstructed images in order to obtain a reliable diagnosis. 

Our results can serve as a basis to develop mathematical criteria to
distinguish added artifacts arising due to the information content in
the data from motion artifacts which occur if the motion is not
correctly compensated for. This can have a great benefit in
applications, for example in the course of estimating the \emph{a
priori} unknown motion parameters which is required in order to apply
a motion compensation algorithm for the reconstruction. To this end,
one first has to develop a motion model which describes the type of
movement performed by the object, and then, the parameters of this
model have to be estimated from the measured data via analytic
\cite{lu_mackie} or iterative \cite{katsevich_local-2011} methods.
However, the estimated parameters will always be affected by errors,
especially in the iterative procedure.  Hence, motion artifacts as
well as added artifacts described in this article will appear in the
reconstructed images. In this case, it is essential to understand and
evaluate whether any given artifact is related to an inaccurate motion
model and incorrect parameters or whether it is inevitable due to
information missing from the data.

\begin{appendix} 

\section{Appendix}

\subsection{The forward operator: proof of Theorem
\ref{theorem_fio}}\label{proof:fio}

Let $f\in \cD(\rtwo)$ and let $\cF$ be the Fourier
transform on $\rtwo$ and    let $\cF_s$ be the
one dimensional Fourier transform in the $s$ variable with the
following normalizations:
\[\cF f(\xi)= \frac{1}{2\pi}\int e^{-ix\cdot \xi} f(x)\,\mathrm{d}x,
\quad \cF_s g(\vp,\tau)= \frac{1}{\sqrt{2\pi}}\int e^{-i\tau s}
g(\vp,s)\,\mathrm{d}s.\] Using the {\it Fourier slice theorem} for the
classical Radon line transform with fixed $\vp$, 
\[ \cF (\cRG f) (\varphi, \sigma) = \cF_s (\cR(f\circ \mf)) (\varphi,\sigma) = \sqrt{2\pi} \, \cF (f \circ \mf) (\sigma \theta(\varphi)).\]
Due to this relation and the substitution $z:=\mf x$, we obtain the
following representation
\begin{align*}
\cRG f (\varphi,s) &= (2\pi)^{-1/2} \, \int_\R e^{i\sigma s} \, \cF_s
(\cRG f) (\varphi,\sigma) \, \mathrm{d}\sigma\\
&= \int_\R e^{i\sigma s} \cF(f\circ \mf) (\sigma \theta(\varphi)) \,
\mathrm{d} \sigma\\
&= (2\pi)^{-1} \, \int_\R e^{\mathrm{i} \sigma s} \int_{\R^2} f(\mf x)
\, e^{- \mathrm{i} \sigma x^T \theta(\varphi)} \, \mathrm{d}x
\mathrm{d} \sigma\\
&= (2\pi)^{-1} \, \int_\R e^{\mi \sigma s} \int_{\R^2} f(x) \, |\det
D\imf x| \, e^{- \mi \sigma (\imf x)^T\thetaphi} \, \md x \md \sigma\\
&= \int_\R \int_{\R^2} e^{\mi \sigma(s-(\imf x)^T\thetaphi)} \, f(x) |\det D \imf x| \, (2\pi)^{-1} \, \md x \md \sigma.
\end{align*}

The function \[\Phi(\varphi,s,x,\sigma)=\sigma(s-(\imf
x)^T\theta(\varphi)) =\sigma(s-H(\varphi,x))\] is homogeneous of degree $1$ with respect to
$\sigma$.  A calculation using this definition shows \begin{align*}
\partial_\sigma \Phi &= \paren{s-(\imf
x)^T\theta(\varphi)}\md \sigma=\paren{s-H(\varphi,x)}\md \sigma,\\
\partial_s \Phi &= \sigma\md s,\\
\partial_x \Phi &= - \sigma \left( (D_x \imf x)^T \right)
\theta(\varphi) \, \md x = -\sigma \cN(\vp,x)\md x\end{align*}
which we justify using \eqref{abbreviation_H} and \eqref{def:N}.  Since
$\Gamma_\varphi$ is a diffeomorphism, the Jacobian matrix
$D_x\paren{\imf x}$ has nowhere zero determinant, so the product
$\paren{D_x\paren{\imf x}}^T\theta(\vp)$ is nowhere zero.  Thus,
altogether, we obtain that $(\partial_{(\varphi,s)} \Phi,
\partial_\sigma \Phi)$ and $(\partial_x \Phi,\partial_\sigma \Phi)$
are nonzero for all $(\varphi,s,x,\sigma)$. Hence, $\Phi$ is a phase
function.  Note that $\Phi$ is nondegenerate because
$\frac{\partial}{\partial s}\paren{\frac{\partial}{\partial
\sigma}\Phi}=1$ is nonzero.

Since $\Gamma_\varphi$ and its inverse are smooth in $(\vp,x)$, the
amplitude of $\cRG$, $a=\abs{D_x\paren{\imf x}}$, and phase function,
$\Phi$, are smooth on their respective domains.  Furthermore,
$a(\varphi,s,x,\sigma)$ does not depend on $\sigma$, so it is
trivially a symbol of order $0$ (see \eqref{def:symbol order}).  This
means that $\cRG$ is an FIO with order $-1/2$. Since $\Gamma_\varphi$
is a diffeomorphism for each $\vp\in \netpe$, the symbol $a$ is
positive and bounded away from zero on every compact set in
$\netpe\times\rtwo$ (and arbitrary $\sigma$).  This shows that the
amplitude $a$ is elliptic and so $\cRG$ is an elliptic FIO.

\subsection {The forward operator: proof of Theorem
\ref{theorem_wavefrontsets}}\label{proof:wavefrontsets}

 According to Theorem \ref{theorem_fio}, $\cRG$ is a Fourier integral
operator. Thus, (\ref{wavefrontset_general_motion}) follows by the
H\"ormander-Sato Lemma \ref{theorem:HS}. 

Now assume the motion model in addition fulfills the Bolker
assumption.  As noted in Theorem \ref{theorem_fio}, the symbol of
$\cRG$ is elliptic.  The proof of the theorem in full generality
follows from the general calculus of FIO in \cite{Ho1971} and it will
be outlined.

Let $f\in\cE'(\rtwo)$ and let $(\xo,\xio)\in \WF(f)\cap
\cV_{\netpe}$.  Then, the set 
\[\Cxo=\Pi_R\inv\sparen{(\xo,\xio)}\]
is nonempty.  By the Bolker Assumption $\Pi_L$ is an immersion 
and so $\Pi_R$ is also an immersion by Prop.\ 4.1.3 \cite{Ho1971}.
Therefore, $\Cxo$ is a discrete set in $\cCG$.  To better understand
this set, we will use the diffeomorphism $c:\netpe\times \rtwo\times
\paren{\rr\smo}\to \cCG$, given in \eqref{def:c}.  
Let  \begin{align*}\lo&=c\paren{\vpo,\xo,\sigo}\\
&=\paren{\vpo,H(\vpo,\xo), \sigo\paren{
-\partial_\vp H(\vpo,\xo)+\ds},\xo,\sigo \cN(\vpo,\xo)}\in
\Cxo.\end{align*}
Note that $\xio= \sigo \cN(\vpo,\xo)$.
Without loss of generality, assume $\sigo>0$.
Let \[\etao=\sigo\paren{-\partial_\vp H(\vpo,\xo)+\ds}\]

We now prove that there is a neighborhood $U$ of $\vpo$ such that
$\lo$ is the only point in $\Cxo$ with $\vp\in U$.  Assume not; then
there must be a sequence $(\vp_j)$ that converges to $ \vpo$ and
another sequence $(\sigma_j)$ in $\rr\smo$ such that
$\Pi_R\paren{c(\vp_j,\xo,\sigma_j)} =(\xo,\xio)$.  However, a
calculation using the definitions of $\Pi_R$ and $c$ shows that
$\sigma_j= \frac{\norm{\xio}}{\norm{D_x H(\vp_j,\xo)}}$.  Therefore
$\sigma_j\to\sigo$ and $c(\vp_j,\xo,\sigma_j)\to
c(\vpo,\xo,\sigo)=\lo$ in $\Cxo$.  This contradicts the fact that
$\Cxo$ is discrete.

Let $\phi_0$ be a smooth cutoff function supported in $U$ and equal to
one in a smaller neighborhood of $\vpo$, and let $\phi_1$ be a cutoff
function equal to one in a neighborhood of $\so = H(\vpo,\xo)$.  For
$(\vp,s)\in \netper$ let $\phi(\vp,s) = \phi_0(\vp)\phi_1(s)$.  Now,
let \bel{def:Mphi}\Mphi (g) = \phi g.\ee Then, $\Mphi:\cD'(\netper)\to
\cE'(\netper)$ is trivially a pseudodifferential operator that has
amplitude $\phi(\vp,s)$ (that is constant in $\eta$) and is nonzero
and hence elliptic at $(\vpo,\so,\etao)$.

Let $\cRGst:\cE'(\netper)\to\cD'(\rtwo)$ be the formal adjoint of
$\cRG:\cD(\rtwo)\to\cE(\netper)$.   Note that in this non-periodic
case, $\cRGst$ is not the
backprojection defined by \eqref{def:Rt} but the
dual operator defined by \eqref{def:R*}.
Furthermore,  $\cRGst$ is an FIO with canonical relation $\cCGt$.

Because $\phi$ has compact support, $\cRGst$, $\Mphi$ and $\cRG$ can be
composed.  Because $\Pi_L$ is an immersion, $\cCG$ and $\cCGt$ are local
canonical graphs, so the composition $\cRGst\Mphi \cRG$ is an FIO
associated to canonical relation
\[\cCGt\circ \cCG \subset \Delta:=\sparen{(x,\xi;x,\xi)\st
(x,\xi)\in T^*(\rtwo)\smo}.\] Since $\cCGt\circ \cCG \subset \Delta$,
$\cRGst\Mphi\cRG$ is a pseudodifferential operator.

The top order symbol of $\cRGst(\Mphi \cRG)$ at $(\xo,\xio)$ is
essentially \bel{symbol}\phi(\vpo,H(\vpo,\xo)) \frac{\abs{\det(D_x
\Gamma_{\vpo} \xo)}^2}{2\pi\norm{\xio}}\ee as can be shown using the
symbol calculation in the proof of Theorem 2.1 in \cite{Q1980}.  Also,
as $\Pi_R:\cCG\to T^*(\rtwo)\smo$ is a conic immersion, the Inverse
Function Theorem shows that $\vp$ is a smooth function of $(x,\xi)$ at
least for $\vp$ near $\vpo$ and for $x$ near $\xo$.  Note that we use
that this symbol is nonzero on only one element of $\Cxo$, $\lo$,
since $\vpo$ is the only angle in $U$ associated to an element of
$\Cxo$.  This symbol is elliptic near $(\xo,\xio)$ because it is
nonzero and homogeneous in $\xi$.  Therefore, $\cRGst(\Mphi \cRG)$ is
elliptic near $(\xo,\xio\dx)$.  So, as $(\xo,\xio\dx)\in
\WF(f)$,\[(\xo,\xio)\in \WF(\cRGst(\Mphi \cRG)).\]

Let $\Pi_L^t:\cCGt\to T^*(\netper)$ and $\Pi_R^t:\cCGt\to T^*(\rtwo)$
be the natural projections.  Since \[(\xo,\xio\dx)\in
\WF\paren{\cRGst\bparen{\Mphi \cRG(f)}}\subset \cCGt\circ
\WF(\Mphi\cRG f)=\Pi_R^t\paren{\paren{\Pi_L^t}\inv\paren{\WF(\Mphi\cRG
f)}},\] some element of $\Pi_L^t(\Cxo)$ is in $\WF(\Mphi \cRG f)$.
Since $\lo$ is the only covector in $\Cxo$ on which the symbol of
$\cRGst M_{\phi}\cRG$ is nonzero, $\Pi_L(\lo)=(\vpo,H(\vpo,\xo),\etao)$
is the only element of $\Pi_L^t(\Cxo)$ on which $\Mphi$ is nonzero.
Therefore, $(\vpo,H(\vpo,\xo),\etao)\in \WF(\cRG f)$.

\subsection{The smoothly periodic case: proof of Theorem
  \ref{thm:psido-R}}\label{proof:psido-R}

The proof of the theorem in full generality follows from arguments in
\cite{Gu1975, GS1977, Q1980}.

Since the motion model is smoothly periodic, we can use Proposition
\ref{prop:composition} to infer $\cRG:\cE'(\rtwo)\to\cE'(\otpr)$ and
$\cRGt:\cD'(\otpr)\to\cD'(\rtwo)$ (which is the formal adjoint in
this case) are both continuous and they can be composed with any
pseudodifferential operator $\cP:\cE'(\otpr)\to \cD'(\otpr)$. 

We first show \bel{surjective} \Pi_R:\cCG\to T^*(\rtwo)\smo \ \ \text{
is surjective.}\ee This will imply that \[
\Pi_R\paren{\Pi_L\inv\paren{T^*(\otpr)\smo}}=T^*(\rtwo)\smo,\] so, from
the discussion in Section \ref{sect:visible singularities},
$\cV_{\otp}=T^*(\rtwo)\smo$ and every singularity is visible.

By \eqref{immersion_condition}, $D_x H(\vp,x)$ is never zero (or the
determinant $\IC(x,\vp)$ would be zero).  For the same reason,
$D_\vp\paren{D_x H(\vp,x)}$ is never zero and $D_x H(\vp,x)$ and
$D_\vp\paren{D_x H(\vp,x)}$ are not parallel.  

Fix $\xo\in \rtwo$.  Consider the function $A:\otp \to S^1$ defined by
\[A(\vp):=\frac{D_x H(\vp,\xo)}{\norm{D_x H(\vp,\xo)}}\in S^1.\] The
map $A$ is periodic of period $2\pi$ and continuous since the motion
model is smoothly periodic.  Because $D_x H(\vp,\xo)$ and $D_\vp
\paren{D_x H(\vp,\xo)}$ are not parallel, a calculus exercise shows
that $A'(\vp)$ is never zero.  Therefore, the $2\pi$ periodic path
\[[0,2\pi]\ni\vp \mapsto A(\vp)\in S^1\] starts at $A(0)$ and ends at
$A(2\pi)=A(0)$ and moves in only one direction.  This shows that the
range of $\vp\mapsto A(\vp)$ is all of $S^1$.  

 Let $\xo\in \rtwo$ and $\xio\in \rtwo\smo$.  Let $\vpo\in\otp$ be an
angle so that $D_x H(\vpo,\xo)$ is parallel to $\xio$.  This can be
done because $\vp\mapsto A(\vp)$ has range $S^1$.  In the global
coordinates on $\cCG$ given by \eqref{def:c},
\bel{PiR}\Pi_R\paren{c(\vpo, \xo, \sigma)}= \paren{\xo,\sigma
\cN(\vpo,\xo)}\ee and for appropriate $\sigma\neq 0$, $\sigma D_x
H(\vpo,\xo)=\xio$.  Therefore $\Pi_R:\cCG\to T^*(\rtwo)\smo$ is
surjective.

Furthermore, because $A'(\vp)$ is never zero and $\otp$ is compact,
there are at most a finite number of angles $\vp\in \otp$ with
$A(\vp)=\xio/\norm{\xio}$.  This shows that there are only a finite
number of points in $\cCG$ that map to $(\xo, \xio)$.  (Here one can
use \eqref{PiR} to show that, for each $(\vp,\xo)$, $\sigma \mapsto
\Pi_R\paren{c(\vp,\xo,\sigma)}$ is one-to-one.) 

Now, we prove the theorem.  Because $\Pi_R$ is surjective and
$\Pi_L$ is injective, 
$\cCGt\circ \cCG = \Delta.$  Because $\cCG$ and $\cCGt$ are local
canonical graphs and $\cRGst$, $\cP$, and $\cRG$ can be composed as
FIO, the composition \[\cL = \cRGst \cP\cRG\] is a pseudodifferential
operator. 

We now explain why $\cL$ is elliptic.  Let $(\xo,\xio)\in
T^*(\rtwo)\smo$.  By the discussion about the map $A$ above, there are
a finite number of angles $\sparen{\vp_0,\dots,\vp_N}$ such that
$\Pi_R\paren{c(\vp_j,\xo,\sigma_j)}=(\xo,\xio)$.  

The symbol of $\cRG$ at $c(\vp_j,\xo,\sigma_j)$ is $a=
\abs{D_x\Gamma_{\vp_j} \xo}$ (see \eqref{amplitude-R}) and the symbol
of $\cRGst$ is the same \cite{Ho1971}.  Let $p$ be the symbol of
$\cP$. Then, by the calculus of FIO, the top order symbol of $\cL$ at
$(\xo,\xio)$ is the sum of $a^2 p/\norm{\xi}$ summed at each element
of the finite set \bel{def:S}S=\sparen{c(\vp_j,\xo,\sigma_j)\st
j=0,\dots,N}.\ee
The proof this statement is completely analogous to the proof of
Theorem 2.1 and equation (15) in \cite{Q1980}.

Since each term in this finite sum is positive as the symbol  $p$ is
everywhere positive and elliptic, the symbol of $\cL$ is
positive.  Therefore, $\cL$ is an elliptic pseudodifferential operator
(the complete argument is analogous to the symbol calculation in the
proof of Theorem 2.1 in \cite{Q1980}).  This proves our theorem.

\begin{Remark}\label{remark:elliptic general}
  Looking over the end of the proof of Theorem \ref{thm:psido-R}, one
sees that the condition for ellipticity is fulfilled as long as the
sum of $a^2 p/\norm{\xi}$ evaluated at each element of the finite set
$S$ given by \eqref{def:S} is an elliptic symbol.

This discussion shows that $\cP$ needs to be elliptic only on
$\Pi_L(\cCG)$, since $S$ is the only set at which the symbol is summed,
and $S$ is a subset of $\cCG$, so its symbol $p$ is only evaluated on
points in $\Pi_L(\cCG)$.  Examples of such pseudodifferential
operators are the operator of Lambda tomography, $-d^2/ds^2$ and the
standard filtered backprojection filter for the linear Radon line
transform, $\sqrt{-d^2/ds^2}$.
\end{Remark}

\subsection{The non-periodic case: Proofs of Theorems
\ref{thm:artifacts} and \ref{thm:ellipticity
nonperiodic}}\label{proof:artifacts}

 \textbf{Proof of Theorem \ref{thm:artifacts}} We apply a
paradigm given in \cite{FrikelQuinto2015} that characterizes the
visible and added singularities in a broad range of incomplete data
tomography problems.  The paradigm uses the following result, which is
a special case of a result of H\"ormander's \cite{Hoermander03}. 

 \begin{lemma} \label{lemma:Ho} Let $u\in \cE'(\netper)$ and let $B$
be a closed subset of $\netper$ with nontrivial interior.  If the
following non-cancellation condition holds \bel{non-cancellation}
\forall (y,\xi)\in \WF(u),\ (y,-\xi)\notin \WF(\chi_B),\ee then the
product $\chi_B u$ can be defined as a distribution.  In this case, 
\[WF(\chi_B u)\subset \cQ(B,\WF(u))\]  where for $W\in T^*(\netper)$
\bel{def:Q}\begin{aligned}\cQ(B,W) :=& \big\{(y,\xi+\eta)\st y\in
B\,, \bparen{(y,\xi)\in W\text{\rm\ or } \xi = 0}\\&\qquad
\text{\rm\ and } \big[(y,\eta)\in \WF(\chi_B)\text{\rm\ or } \eta =
0\big]\big\}\,.\end{aligned}\ee
\end{lemma}

To prove Theorem \ref{thm:artifacts}, we apply this paradigm to the
Fourier integral operator $\cRGc$ with the data set $B:= [0,2\pi]
\times \R$.  We first use this lemma to establish that the operator
$\cLr$ is well defined.

 \begin{proposition}\label{prop:composition nonpdic} For $f\in \cE'(\rtwo)$,
$\cotpr$ can be multiplied by $\cRG f$ as distributions.  Let $\psi$
be a smooth function equal to $1$ on $\otp$ and supported in $\netpe$
and let $\cRGtr=\cRGst\psi$.  Then, for $\cP$ a pseudodifferential
operator, $\cRGtr$, $\cP$ and $\cotpr\cRG$ can all be composed and
$\cLr$ given in \eqref{def:cLr} is defined and $\cLr:\cE'(\rtwo)\to
\cD'(\rtwo)$.\end{proposition}

\textbf{Proof:} First, we show that $\cP \cRGr f$ is a distribution.
The product $\cotpr \cRG f$ is well-defined for distributions $f\in
\cE'(\R^2)$, since $\WF(\cotpr)$ has $\md s$ component of zero,
whereas any covector in $\cCG \circ \WF(f)$ has nonzero $\md s$
component by the definition of $\cCG$, \eqref{def:C}.  Therefore, the
non-cancellation condition in Lemma \ref{lemma:Ho} holds and $\cotpr
\cRG f$ is a distribution.  

We claim $\cotpr\cRG f$ has compact support.  First, this distribution
has support in $\otpr$ because $\cotpr$ does.  Since, for each $\vp$,
$s\mapsto C(\vp,s)$ is a smooth foliation of the plane, for each
$\vp$, the support in $s$ of $\cotpr\cRG f(\vp,\cdot)$ is compact.
Since the foliation depends smoothly on $\vp$ and $\vp$ is in the
compact set $\otp$, there is an $M>0$ such that the support of
$\cotpr\cRG f$ is in $\otp\times [-M,M]$.  Therefore,  $\cP
\cRGr f$ is defined as a distribution in $\cD'(\netper)$.

One proves that $\psi \cRG$ is continuous from $\cD(\rtwo)$ to
$\cD(\netper)$ using the same arguments as in the proof of Proposition
\ref{prop:composition}.  This implies that
$(\psi\cRG)^* = \cRGst \psi=\cRGtr$ is weakly continuous from
$\cD'(\netper)$ to $\cD'(\rtwo)$.  Therefore, $\cLr f$ is defined as a
distribution.\qed

 We continue the proof of Theorem \ref{thm:artifacts} and now use
Theorem \ref{theorem:HS} to show \bel{composition1} \WF(\cRG
f)\subset \cCG \circ \WF(f).\ee Next, we use Lemma \ref{lemma:Ho} to
get an upper bound for $\WF(\cP\cRGr f)$.  Using \eqref{def:Q} and
\eqref{composition1}, we obtain 
\[\WF(\cP \cRGr f)\subset \cQ\paren{\otpr,\WF(\cRG f)}\subset
\cQ\paren{\otpr, \cCG\circ \WF(f)}.\] Then, \begin{align*} \cQ (\otpr,
\cCG \circ \WF(f)) = \big[\paren{\cCG \circ \WF(f)}
\cap&T^*_{\otpr}(\netper) \big] \\ & \cup \WF(\cotpr) \cup W_{\lbrace
0,2\pi\rbrace} (f),
\end{align*}
where $T^*_{\otpr}(\netper)$ is defined in \eqref{def:T*B}
\begin{align*} W_{\lbrace 0,2\pi\rbrace} (f) = \Big\lbrace
(\varphi,s,\sigma \md s & + [\mu - \sigma \partial_\varphi H(\vp,x)]\md
\varphi)\st  \\ 
& \sigma, \mu \neq 0, \varphi \in \lbrace 0,2\pi \rbrace, s \in \R \\
&\qquad  x \in C(\varphi,s), \ \text{and} \  (x,\sigma
\cN(\vp,x) )\in \WF(f)\Big\rbrace. \end{align*} Equivalently, this set
can be written as
\begin{align} \label{pre_set_add_art} W_{\lbrace 0,2\pi\rbrace} (f) =
\Big\lbrace (\varphi,s,\sigma \md s + \nu \md \varphi)\st \sigma \neq
0,\ & \nu \in \R, \varphi \in \lbrace 0,2\pi \rbrace, s \in \R \\ &
\exists x \in C(\varphi,s), (x,\sigma \cN(\vp,x)) \in
\WF(f)\Big\rbrace. \notag \end{align} 

To accomplish the final step of the paradigm, we determine
\[\cCGt \circ \cQ \paren{\otpr, \cCG \circ \WF(f)},\] which
corresponds to computing the three components
\begin{align*} \cCGt \circ \cQ (\otpr, \cCG \circ \WF(f)) = \cCGt \circ &
                                                                      [(\cCG
\circ \WF(f)) \cap T^*_{\otpr}(\netper) ]\\ & \cup \cCGt \circ
\WF(\chi_A) \\ & \quad \cup \cCGt \circ W_{\lbrace 0,2\pi\rbrace}(f)
.\end{align*} Since $\cCG$ fulfills the Bolker assumption, $\cCGt \circ
\cCG \circ \WF(f) \subset \WF(f)$.  Thus, for the first component, we
obtain
\[ \cCGt \circ  \bparen{(\cCG \circ \WF(f)) \cap T^*_{\otpr}(\netper)} \subset \WF(f) \cap \cV_{[0,2\pi]},\]
i.e. the set of visible singularities, $\WF_{[0,2\pi]}(f)$.

For the second component, $\cCGt \circ \WF(\chi_A) = \emptyset$, since
the $\md s$ component of any covector in $\WF(\chi_A)$ is zero and all
covectors in $\cCGt$ have nonzero $\md s$ component.

Lastly, we consider $\cCGt \circ W_{\lbrace 0,2\pi\rbrace}(f)$ and
show that this equals the set of additional artifacts $\cA(f)$.  To
this end, we let
\[ \rho = (\varphi,s,\nu\md \varphi + \sigma \md s)\in W_{\lbrace 0,2\pi\rbrace}(f),\]
and so $\varphi \in \lbrace 0,2\pi\rbrace, \, s,\nu \in \R, \, \sigma
\neq 0$ and there is $x \in C(\varphi,s)$ such that $(x,\sigma
\cN(\vp,x)) \in \WF(f)$.  Using the definition of composition, one
sees
\[ \cCGt \circ \lbrace \rho \rbrace = \left\lbrace (\tilde{x},
\sigma (\cN(\vp,\tx))) \st (\tilde{x},\sigma \cN(\vp,\tx),\rho) \in
\cCGt \right\rbrace.\] By definition of $\cCGt$, $\tilde{x} \in
C(\varphi,s)$, i.e.  $s=H(\tilde{x},\varphi)$ and $-\nu/\sigma =
D_\varphi H(\tilde{x},\varphi)$. Since $\nu$ is arbitrary, for any
$\tilde{x}$ in $C(\varphi,s)$ there is a corresponding covector in
this composition. Therefore, for any $\tilde{x} \in C(\varphi,s)$, the
covector $(\tilde{x},\sigma \cN(\vp,\tx)) \in \cCGt \circ W_{\lbrace
0,2\pi\rbrace}(f)$. Thus, this set corresponds to the set of possible
added singularities (\ref{set_add_sing}). \qed

\textbf{Proof of Thm.\ \ref{thm:ellipticity nonperiodic}} 
Let $(\xo,\xio)\in \cV_{(0,2\pi)}$, then by the uniqueness assumption
\eqref{unique covector}, there is a unique $(\vpo,\so)\in \netper$
such that $\xio$ is conormal to $C(\vpo,\xo)$ at $\xo$.  Since $\vpo$
is unique and $(\xo,\xio)\in \cV_{(0,2\pi)}$, $\vpo\in (0,2\pi)$. Let
$\sigo$ be the unique nonzero number such that
$\xio=\sigo\cN(\vpo,\xo)$.  Then, \bel{def:lo}\lo
=c(\vpo,\xo,\sigo)\in\cCG\ee is the unique covector in $\cCG$ such
that $\Pi_R(\lo) = (\xo,\xio)$ (where $c$ is given by \eqref{def:c}).
Let \bel{def:rhoo}\rhoo: =\Pi_L(\lo)=
(\vpo,\so,\sigo\paren{-\partial_\vp H(\vpo,\xo) +\md s}).\ee We note
that \bel{unique2}\sparen{\rhoo}= \cCG\circ\sparen{(\xo,\xio)},\quad
\sparen{(\xo,\xio)}=\cCGt\circ \sparen{\rhoo}.\ee These equalities
are true by \eqref{circ and projections} and the Bolker Assumption
because $\lo$ is the only element in $\Pi_R\inv\sparen{(\xo,\xio)}$.

First, we show $\WF_{(0,2\pi)}(\cLr f)\subset \WF_{(0,2\pi)} (f) $.
Assume the covector \[(\xo,\xio)\in \WF_{(0,2\pi)}(\cLr f).\] Using
the result of the last paragraph, let $\vpo\in(0,2\pi)$ and $\sigo\neq
0$ be the unique numbers so that $\xio=\sigo\cN(\vpo,\xo)$.  By
Theorem \ref{thm:artifacts}, in particular \eqref{set_add_sing},
\[(\xo,\xio)\in \WF_{[0,2\pi]}(f)\cup \cA(f).\]  However, $\cA(f)$
includes singularities $(x,\sigma \cN(\vp,x))$ only for $\vp=0$ or
$\vp=2\pi$ and by the uniqueness assumption, \eqref{unique covector},
since $\xio=\sigo\cN(\vpo,\xo)$ and $\vpo\notin \sparen{0,2\pi}$,
$(\xo,\xio)\notin \cA(f)$, so $(\xo,\xio)\in \WF_{(0,2\pi)}(f)$.

Now, let $(\xo,\xio)\in \WF_{(0,2\pi)}(f)$. Ellipticity and the
uniqueness assumption will be used to show that $(\xo,\xio)\in
\WF_{(0,2\pi)}\paren{\cLr(f)}$.  Let $\vpo$, $\so$, $\sigo$, $\lo$,
and $\rhoo$ be as in the first paragraph of this proof for
$(\xo,\xio)$.  As noted above, $\vpo\in (0,2\pi)$ by the uniqueness
assumption.  Let $\Mphi$ be the cutoff operator given by
\eqref{def:Mphi} in the proof of Theorem \ref{theorem_wavefrontsets}.
The function $\phi$ in the definition of $\Mphi$ is the product of two
compactly supported cutoff functions, $\phi_0(\vp)$ and $\phi_1(s)$,
and we assume that the cutoff function at $\vpo$, $\phi_0$, is also
supported in $(0,2\pi)$.  As in the proof of Theorem
\ref{theorem_wavefrontsets}, 
\[\cRGtr \cP\Mphi \cRGr =\cRGst \paren{\psi\cP \Mphi\cotpr \cRG}\] is an elliptic pseudodifferential operator
near $(\xo,\xio)$ because its symbol is
\bel{symbol2}\phi(\vpo,H(\vpo,\xo))p(\rhoo)\frac{\abs{\det(D_x
\Gamma_{\vpo} \xo)}^2}{2\pi\norm{\xio}}\ee where $p$ is the top order
symbol of $\cP$.  (Note that $\Mphi\cotpr = \Mphi$ since the support
of $\phi$ is in $(0,2\pi)\times \rr$.  Also, the cutoff $\psi$ has no
effect on the top order symbol \eqref{symbol2} since $\phi\cdot \psi
=\phi$ as $\psi$ is equal to one in $\otp$.) So
\bel{fact1}(\xo,\xio)\in \WF_{(0,2\pi)}(\cRGtr \cP\Mphi \cRGr f).\ee

We now show that \bel{fact2}(\xo,\xio)\notin
\WF_{(0,2\pi)}\paren{\cRGtr\paren{ \cP M_{(1-\phi)} \cotpr \cRG
f}} \ee by showing 
 \bel{simpler}(\xo,\xio)\notin \cCGt\circ \WF\paren{
\cP M_{(1-\phi)} \cotpr \cRG f} \ee and then using the H\"ormander-Sato
Lemma \ref{theorem:HS}.

Because $(1-\phi)$ is zero near $\vpo$, $M_{(1-\phi)}\cRGr f$ is
microlocally smooth near $\rhoo$.  So, $\psi\cP M_{(1-\phi)}\cRGr f$
is microlocally smooth near $\rhoo$.  But, by \eqref{unique2}, $\rhoo$
is the only covector in $\Pi_L(\cCG)$ that could map to $(\xo,\xio)$
under $\Pi_R^t\circ \Pi_L^t$.  Therefore, \eqref{simpler} holds and this
proves \eqref{fact2}.

Putting \eqref{fact1} and \eqref{fact2} together, we see that
$(\xo,\xio)\in \WF(\cLr f)$, and this finishes the proof.  \qed

\subsection{Our theorems for arbitrary smooth
weights}\label{sect:arbitrary measure} 

Finally, we explain why our theorems are true even if the weight $
|\det D\Gamma_\varphi^{-1} x| $ in the definition of $\cRG$,
\eqref{def:R}, and the definition of $\cRGt$, \eqref{def:Rt}, are
replaced by smooth positive weights.  Basically, this is true because
elliptic FIO associated to the same canonical relation have the same
microlocal properties, and Radon transforms that integrate over the
same sets (associated to the same double fibration \cite[Definition
1.1]{Q1980}) are FIO with the same canonical relations.

Let $\mu$ be a smooth positive function on $\netpe\times \rtwo$, then
\[ {\cRG}_\mu f(\vp,s)=\int_{x\in C(\vp,s)} f(x) \mu(\vp,x)\md x\] is
an elliptic FIO associated to $\cCG$.  This is true by the general
theory of Radon transforms as FIO \cite{Gu1975, GS1977} (see also
\cite{Q1980}) because this transform integrates over the same sets,
$C(\vp,s)$, as $\cRG$ and the weight is smooth and nowhere zero.

In the smoothly periodic case, the weight, $\mu$ for ${\cRG}_\mu$ must
be $2\pi$-periodic.  In this case, a generalized backprojection can
be defined as 
\[\cR^\dagger_{\Gamma\,\nu} g(x) = \int_{\vp\in\otp} 
g(\vp,H(\vp,x))\nu(\vp,x)\md \vp. \] where $\nu$ is a positive
smooth $2\pi$-periodic function.  Because the weights are smooth and
positive $\cRG_\mu$ and $\cR^\dagger_{\Gamma\,\nu}$ are elliptic and
associated to $\cCG$ and $\cCGt$ respectively.  The proof of
Proposition \ref{prop:composition} for $\cR^\dagger_{\Gamma\,
\nu}\cP\cRG_\mu$ does not change, and the other proofs for the
smoothly periodic case rest on the fact that these transforms are
elliptic and associated with the same canonical relations as $\cRG$
and $\cRGt$.

For the non-periodic case, the weighted 
backprojection operator is
\[\cRGtr=\int_{\vp\in\netpe}
\phi(\vp)\nu(\vp,x)g(\vp,H(\vp,x))\md \vp\]
where $\phi$ is a smooth function equal to one on $\otp$ and supported in
$(-\eps,2\pi+\eps)$.  In this case, too, the proofs are the same because
the transforms have the same microlocal properties.

\end{appendix}

\subsection*{Acknowledgments}\label{sect:acknowledgments} The work of
the first author was partially supported by a fellowship of the German
Academic Exchange Service (DAAD) under 91536377. She thanks Tufts
University for its hospitality during her stay as visiting scholar.
The second author thanks Universit\"at des Saarlandes for its
hospitality during multiple visits.  The second author is indebted to
Jan Boman and J\"urgen Frikel for stimulating conversations related to
microlocal analysis and tomography and to Alexander Katsevich for
insightful comments about this article and the general problem. 

The work of the second author was partially supported by NSF grant DMS
1311558.

\end{document}